\documentclass[a4paper,oneside,11pt]{article}

\usepackage{amsfonts,color}

\usepackage[final]{pdfpages}
\usepackage{graphicx, soul}
\usepackage{caption}

\newtheorem{dfn}{Definition}[section]
\newtheorem{tw}[dfn]{Theorem}
\newtheorem{prop}[dfn]{Proposition}
\newtheorem{rem}[dfn]{Remark}
\newtheorem{ex}[dfn]{Example}

\newtheorem{cor}[dfn]{Corollary}
\usepackage{amssymb} 
\usepackage{amsmath}
\numberwithin{equation}{section}

\renewcommand{\theequation}{\thesection.\arabic{equation}}

 \global\long\def\sbr#1{\left[ #1\right] }

 \global\long\def\rbr#1{\left(#1\right)}

 \global\long\def\R{\mathbb{R}}

 \global\long\def\dd#1{\textnormal{d}#1}

 \global\long\def\ra{\rightarrow}
 
 \global\long\def\ns{\infty}

\usepackage{anysize}
\usepackage{array}
\usepackage{enumerate}

\author{Micha\l \ Barski  \\ \small  Faculty of Mathematics, Warsaw University, Poland\\
 \small{\it m.barski@mimuw.edu.pl} \bigskip \\
\\
Rafa\l \ \L ochowski
\\ \small 
Department of Mathematics and Mathematical Economics,\\ \small Warsaw School of Economics, Poland\\ \small{\it rlocho@sgh.waw.pl}}

\title{\bf Classification and calibration of affine models driven by independent L\'evy processes}

\begin{document}

\maketitle

\begin{abstract}
\renewcommand{\theequation}{\arabic{equation}}

The paper is devoted to the study of the short rate equation of the form
\begin{gather}\label{equation abstract}
\dd R(t)=F(R(t))\dd t+\sum_{i=1}^{d}G_i(R(t-))\dd Z_i(t), \quad R(0)=x\geq 0,\quad t>0,
\end{gather}
with deterministic functions $F,G_1,...,G_d$ and independent L\'evy processes of infinite variation  $Z_1,...,Z_d$ with regularly varying Laplace exponents. The equation is supposed 
to have a nonnegative solution which generates an affine term structure model. A  precise form of the generator of $R$ is characterized and
a related classification of equations which generate affine models introduced in the spirit of Dai and Singleton \cite{DaiSingleton}.
Each class is shown to have its own canonical representation which is an  equation with the same drift and the jump diffusion part based on a L\'evy process taking values in $\mathbb{R}^{g}, 1\leq g\leq d$, with independent coordinates being stable processes with stability indices in the range $(1,2]$. Numerical calibration results of canonical representations to the market term structure of interest rates  are presented and compared with the classical CIR model. The paper generalizes the classical results on the CIR model from \cite{CIR}, as well as on its extended version from \cite{BarskiZabczykCIR} and \cite{BarskiZabczyk} where $Z$ was a one-dimensional L\'evy process. 
\end{abstract}

\tableofcontents

\section{Introduction}

The study of continuous state branching processes with immigration (CBI) by Kawazu and Watanabe \cite{KawazuWatanabe}
revealed attractive analytical properties of affine processes which motivated  
Filipovi\'c to bring them, in the pioneering paper \cite{FilipovicATS}, in the field of finance.  Affine processes are
widely used in various areas of mathematical finance. They appear in term structure models, by credit risk modelling and are applied within the stochastic volatility framework. Solid fundamentals of affine processes in finance were laid down by Filipovi\'c \cite{FilipovicATS} and by Duffie, Filipovi\'c and Schachermeyer \cite{DuffieFilipovicSchachermeyer}. 
The results obtained in these papers settled a reference point for further research and proved the usefulness and strength  
of the Markovian approach. Missing questions on regularity and existence of c\`adl\`ag versions were answered by Cuchiero, Filipovi\'c and Teichmann \cite{CuchieroFilipovicTeichmann} and  Cuchiero and Teichmann \cite{CuchieroTeichmann}.

The systematic study of affine processes in finance was motivated by classical 
stochastic short rate models, like CIR (Cox, Ingersoll, Ross) \cite{CIR}, Vasi\v cek \cite{Vasicek} and model with diffusion factors of Dai and Singleton \cite{DaiSingleton}, and resulted in discovering new stochastic equations, also with jumps; see, among others, \cite{FilipovicATS}, Duffie and G\^arleanu \cite{DuffieGarleanu},  Barndorff-Nielsen and Shephard \cite{Bandorff-NielsenShepard}, Jiao, Ma and Scotti \cite{JiaoMaScotti}. Nevertheless, the full description of affine processes representable in terms of stochastic equations is far from being clear. This is because the Markovian description of affine processes based on  generators does not, in general, allow encoding the form of a possible underlying stochastic equation. 
The framework based on stochastic dynamics offers, however, unquestionable advantages like discretization schemes enabling 
Monte Carlo simulations which are essential for example for pricing exotic, i.e. path-dependent, derivatives. A comprehensive treatment of simulating schemes for affine processes and pricing methods can be found in \cite{Alfonsi1}. Stochastic equations allow also identifying the number of random sources in the model which is of some use by calibration and hedging. In this paper we focus on recovering from the Markovian setting those affine processes which are given by stochastic equations driven by a multidimensional L\'evy process with independent coordinates. Specifically, we focus on the equation 
\begin{gather}\label{rownanie 1}
\dd R(t)=F(R(t))\dd t+\sum_{i=1}^{d}G_i(R(t-))\dd Z_i(t), \quad R(0)=x,\quad t>0,
\end{gather}
where $x$ is a nonnegative constant, $F$, $\{G_i\}_{i=1,2,...,d}$ are deterministic functions and $\{Z_i\}_{i=1,2,...,d}$ are independent L\'evy processes and martingales. A solution $R(t), t\geq 0$, if nonnegative, will be identified here with the short rate process which defines the bank account process by
$$
B(t):=e^{\int_{0}^{t}R(s)ds}, \quad t\geq 0.
$$
Related to the savings account are zero coupon bonds. Their prices form a family of stochastic processes
$P(t,T), t\in [0,T]$, parametrized by their maturity times $T\geq 0$.  The price of a bond with maturity $T$ at time $T$ is equal to its
nominal value, typically assumed, also here, to be $1$, that is $P(T,T)=1$. The family of bond prices is supposed to have the {\it affine structure}, which means that 
\begin{gather}\label{affine model}
P(t,T)=e^{-A(T-t)-B(T-t) R(t)}, \quad 0\leq t\leq T,
\end{gather}
for some smooth deterministic functions $A$, $B:[0,+\infty) \rightarrow \mathbb{R}$. Hence, the only source of randomness in the affine model \eqref{affine model} is the short rate process $R$ given by \eqref{rownanie 1}.  As the resulting market constituted by $(B(t), \{P(t,T)\}_{T\geq 0})$
should exclude arbitrage, the discounted bond prices
$$
\hat{P}(t,T):=B^{-1}(t)P(t,T)=e^{-\int_{0}^{t}R(s)ds-A(T-t)-B(T-t)R(t)}, \quad 0\leq t\leq T,
$$
are supposed to be local martingales for each $T\geq 0$. This requirement affects in fact our starting equation \eqref{rownanie 1}. Thus the functions $F$, $\{G_i\}_{i=1,...,d}$ and the noise $Z=(Z_1,...,Z_d)$ should be chosen such that $\eqref{rownanie 1}$ has a nonnegative solution for any $x\geq 0$ and such that, for some functions $A$, $B:[0,+\infty) \rightarrow \mathbb{R}$ and each $T\geq 0$, \ $\hat{P}(t,T)$ is a local martingale on $[0,T]$. If this is the case, \eqref{rownanie 1} will be called to {\it generate an affine model} or to be a {\it generating equation}, for short.

The description of all generating equations with one-dimensional noise is well known, see Section \ref{Low-dimensional generating equations} for a brief summary. This paper deals with \eqref{rownanie 1} in the case $d>1$. The multidimensional setting makes the description of generating equations more involved due to the fact that two apparently different generating equations may have solutions which are Markov processes with identical generators. For brevity, we will call such solutions 'identical' or 'the same solutions'. The resulting bond markets are then the same, so such equations can be viewed as equivalent. This phenomenon does not appear in the one-dimensional case, but was a central point in the study of a multi-factor affine models by Dai and Singleton \cite{DaiSingleton}. Recall, in the class of affine models considered in \cite{DaiSingleton} the short rate is an affine function of $N$ factors $(Y_1,...,Y_N):=Y$, which are given by a diffusion equation of the form
\begin{gather}\label{rownanie Dai Singleton}
\dd Y(t)=H(Y(t))\dd t+\Sigma \sqrt{\text{diag} (A+BY(t))}\dd W(t), 
\end{gather}
where $H$ is a specific affine function, $\Sigma, B$ are $N\times N$ matrices, $A$ is a vector in $\mathbb{R}^N$ and the value of $\text{diag}(v)$ is the diagonal $N\times N$ matrix with the coordinates of $v\in \mathbb{R}^N$ on the diagonal. Above $W$ stands for the Wiener process in $\mathbb{R}^N$. By particular choices of parameters, one may recognize in \eqref{rownanie Dai Singleton} many specific models used in practice, for details see \cite{DaiSingleton}. The question of characterization of 
equations \eqref{rownanie Dai Singleton} which generate affine models was handled in \cite{DaiSingleton}, see also \cite{CheriditoFilipovicKimmel}, by classifying the structure of factors. The classification is based on the parameter $m:=\text{rank}(B)$ interpreted as a degree of dependence of the conditional variances on the number of factors. Each equation \eqref{rownanie Dai Singleton} which generates an affine model is classified as a member of one of $N+1$ disjoint subfamilies 
$$
\mathbb{A}_m(N),\quad  m=0,1,...,N,
$$ 
of equations. All equations within a chosen subfamily provide the same short rate and the short rates differ across subfamilies. 
Moreover, each subfamily is shown to have its own {\it canonical representation} for which \eqref{rownanie Dai Singleton} simplifies, i.e. the diffusion matrix in  \eqref{rownanie Dai Singleton} is diagonal.  Although our setting based on equation \eqref{rownanie 1} differs, our approach of characterizing generating equations has much in common with that of Dai and Singleton. The main results of the paper, i.e.
Theorem \ref{TwNiez}, Corollary \ref{cor o postaci mu} and Proposition \ref{prop canonical representation} imply that under mild assumptions any generating equation \eqref{rownanie 1} has the same solution as that of the following equation 
\begin{gather}\label{canonical representation rownania}
\dd R(t)=(a R(t)+b)\dd t+\sum_{k=1}^{g} d_k^{1/\alpha_k} R(t-)^{1/\alpha_k} \dd Z^{\alpha_k}_k(t),
\end{gather}
with some $1\leq g\leq d$ and parameters  $a\in\mathbb{R}$, $b\geq 0$, $d_k>0$, $k=1,2,...,g$,  driven by independent
stable processes  $\{Z_k^{\alpha_k}\}$ with indices $\{\alpha_k\}$ such that $2\geq \alpha_1>\alpha_2>...>\alpha_g>1$.  All generating equations having the same solutions as 
\eqref{canonical representation rownania} form a set which we denote by
\begin{gather}\label{podklasy dla rownann}
\mathbb{A}_g(a,b;\alpha_1,\alpha_2,...,\alpha_g; \eta_1,...,\eta_g),
\end{gather}
where $\eta_i:=\frac{\Gamma(2-\alpha_i)}{\alpha_i(\alpha_i-1)} d_i, i=1,...,g$, $\Gamma(\cdot)$ is the Gamma function. We call \eqref{canonical representation rownania} a {\it canonical representation} of \eqref{podklasy dla rownann}.
By changing values of the parameters in \eqref{podklasy dla rownann} one can thus
split all generating equations into disjoint subfamilies with a tractable canonical representation for each of them. 

The number and structure of generating equations which form \eqref{podklasy dla rownann} depend on the noise dimension in \eqref{rownanie 1}. As one may expect, the set \eqref{podklasy dla rownann} is getting larger as $d$ increases.
In Section \ref{section Generalized CIR equations on a plane} we determine all generating equations on a plane by formulating
concrete conditions for $F, G$ and $Z_1,Z_2$ in \eqref{rownanie 1}. For $d=2$ the class $\mathbb{A}_1(a,b;\alpha_1; \eta_1)$ 
consists of a wide variety of generating equations while $\mathbb{A}_2(a,b;\alpha_1,\alpha_2; \eta_1,\eta_2)$ turns out to be a singleton.
The passage to the case $d=3$ makes, however, $\mathbb{A}_2(a,b;\alpha_1,\alpha_2; \eta_1,\eta_2)$
a non-singleton. This phenomenon is discussed in Section \ref{section Example in higher dimensions}.

A tractable form of canonical representations is supposed to be an advantage for applications. 
One finds in \eqref{canonical representation rownania} with $g=1,\alpha_1=2$ the classical CIR equation 
and may expect that additional stable noise components improve the model of bond market. 
For $g=2, \alpha_1=2$ and $1<\alpha_2\leq 2$ equation \eqref{canonical representation rownania} 
becomes the alpha-CIR equation studied in \cite{JiaoMaScotti}. It was shown in \cite{JiaoMaScotti} that empirical behaviour of the European sovereign bond market is closer to that implied by the alpha-CIR equation than by the CIR equation
due to the permanent overestimation of the short rates by the latter one. The alpha-CIR equation allows also reconciling low interest rates with large fluctuations related to the presence of jump part whose tail fatness is controlled by the parameter $\alpha_2$.
In the last part of the paper we focus on the calibration of canonical representations to market data. Into account are taken 
the spot rates of European Central Bank implied by the $AAA$ - ranked bonds, Libor rates and six-month swap rates. 
We compute numerically the fitting error for \eqref{rownanie 1}  in the Python programming language with $g$ in the range from $1$ up to $5$. 
This illustrates, in particular, the influence of $g$ on the reduction of fitting error which is always less than in the CIR model. The freedom 
of choice of stability indices makes the canonical model curves more flexible, hence with shapes better adjusted to the market curves.
The effect is especially visible for market data after March 2022 when the curves started to change their shapes.

The structure of the paper is as follows. In Section \ref{section Preliminaries} we discuss 
the Laplace exponents of L\'evy processes, in particular, the Laplace exponents of  the projections of $Z$ along $G$, defined as the processes 
\begin{gather}\label{noise projection introduction}
Z^{G(x)}(t):=\sum_{i=1}^{d}G_i(x)Z_i(t) , \quad t\geq0, \quad x\geq 0,
\end{gather}
which play a central role in the sequel. The second part of Section \ref{section Preliminaries} is based on the preliminary characterization of generating equations, i.e. Proposition \ref{prop wstepny}, which is a version of the result from  \cite{FilipovicATS} characterizing the generator of a Markovian short rate. This leads to a precise formulation of the problem studied in the paper. Further we 
describe one dimensional generating equations and discuss the non-uniqueness of generating equations in the multidimensional case.
In Example \ref{ex different eq ident sol} we show two different equations with the same solutions. Section \ref{section Classification of generating equations} is concerned with the classification of generating equations. Section \ref{section Noise with independent coordinates} contains the main results of the paper which  provide a precise description of the generator of \eqref{rownanie 1}. This makes more specific the, rather abstract, result from \cite{FilipovicATS} and motivates introducing the classification of generating equations. The required assumption on the Laplace exponent of the noise to vary regularly at zero is reformulated in terms of L\'evy measure in Section \ref{sec slowly varying Laplace exp.}. Section \ref{section Generalized CIR equations on a plane} and Section \ref{section Example in higher dimensions} are devoted to generating equations on a plane and an example in the three-dimensional case, respectively. In Section \ref{section Applications} we discuss the calibration of canonical representations.

\section{Preliminaries}\label{section Preliminaries}
In this section we recall some facts on L\'evy processes needed in the sequel and present a version of the result on generators of
Markovian affine processes \cite{FilipovicATS}, see Proposition \ref{prop wstepny}, which is used for a precise formulation of the problem considered in the paper. We explain the meaning of the projections of the noise \eqref{noise projection introduction}
and show in Example \ref{ex different eq ident sol} two different generating equations having the same projections, hence identical solutions.
For illustrative purposes we keep referring to the one-dimensional case where the forms of generating equations are well known, see Section
\ref{Low-dimensional generating equations} below. For the sake of notational convenience we often use a scalar product notation $\langle\cdot,\cdot\rangle$ in $\mathbb{R}^d$ and write \eqref{rownanie 1} in the form
 \begin{gather}\label{rownanie 2}
\dd R(t)=F(R(t))\dd t+\langle G(R(t-)),\dd Z(t)\rangle, \quad R(0)=x\geq 0, \qquad t>0,
\end{gather}
where $G:=(G_1,G_2,...,G_d):[0,+\infty)\longrightarrow\mathbb{R}^d$ and 
 $Z:=(Z_1,Z_2,...,Z_d)$  is a L\'evy process in $\mathbb{R}^d$.

\subsection{Laplace exponents of L\'evy processes}
Let $Z$ be an $\mathbb{R}^d$-valued L\'evy process  with characteristic triplet $(a,Q,\nu(\dd y))$. Recall, $a\in\mathbb{R}^d$ describes the drift part of $Z$, $Q$ is a non-negative, symmetric, $d\times d$ covariance matrix, characterizing the coordinates' covariance of the Wiener part $W$ of $Z$, and $\nu(\dd y)$ is a measure on $\mathbb{R}^d\setminus\{0\}$ describing the jumps of $Z$. It is called the L\'evy measure of $Z$ and satisfies the condition
\begin{gather}\label{warunek na miare Levyego 1}
\int_{\mathbb{R}^d}(\mid y\mid^2\wedge \ 1)\ \nu(\dd y)<+\infty.
\end{gather} 
Recall, $Z$ admits a representation as a sum of four independent processes of the form
\begin{gather}\label{LevyIto}
Z(t)=at +W(t)+\int_{0}^{t}\int_{\{\mid y\mid\leq 1\}}y\tilde{\pi}(\dd s,\dd y)+\int_{0}^{t}\int_{\{\mid y\mid> 1\}}y\pi(\dd s,\dd y),
\end{gather}
called the L\'evy-It\^o decomposition of $Z$. Above $\pi(\dd s,\dd y)$ and $\tilde{\pi}(\dd s,\dd y):=\pi(\dd s,\dd y)-\dd s \nu(\dd y)$ stand for the jump measure and the compensated jump measure of $Z$, respectively. If
\begin{gather}\label{wariacja Z}
\int_{\{\mid y\mid<1\}}\mid y\mid \nu(\dd y)=+\infty,
\end{gather}
then $Z$ is of infinite variation. If \eqref{wariacja Z} does not hold and $Z$ has no Wiener part, the variation of $Z$ is finite. The coordinates 
of $Z$ are independent if and only if $Q$ is diagonal and $\nu(\dd y)$ is concentrated on axes.

We consider the case when $Z$ is a martingale and call it a L\'evy martingale for short. Its drift and the L\'evy measure are such that 
\begin{gather}\label{warunek na miare Levyego 2}
\int_{\{\mid y\mid>1\}}\mid y\mid\ \nu(\dd y)<+\infty, \quad a+\int_{\{\mid y\mid>1\}}y \ \nu(\dd y)=0.
\end{gather}
Consequently, the characteristic triplet of $Z$ is 
\begin{gather}\label{chrakterystyki Z}
\left(-\int_{\{\mid y\mid>1\}}y \ \nu(\dd y), \ Q, \ \nu(\dd y)\right),
\end{gather}
and \eqref{LevyIto} takes the form
$$
Z(t)=W(t)+X(t), \qquad X(t):=\int_{0}^{t}\int_{\mathbb{R}^d}y \ \tilde{\pi}(\dd s,\dd y), \quad t\geq 0,
$$
where $W$ and $X$ are independent. The martingale $X$ will be called the jump part of $Z$. Its Laplace exponent $J_{X}$,  defined by the equality 
\begin{equation} 
\mathbb{E}\sbr{e^{-\langle \lambda,X(t)\rangle}}=e^{tJ_{X}(\lambda)}, 
\end{equation}
has the  following representation
\begin{equation} \label{Jdef}
J_X(\lambda)=\int_{\mathbb{R}^d}(e^{-\langle\lambda,y\rangle}-1+\langle\lambda,y\rangle)\nu(\dd y),
\end{equation}
and is finite for $\lambda\in\mathbb{R}^d$ satisfying
$$
\int_{\mid y\mid>1}e^{-\langle \lambda,y\rangle}\nu(\dd y)<+\infty.
$$
By the independence of $X$ and $W$ we see that
$$
\mathbb{E}\sbr{e^{-\langle \lambda,Z(t)\rangle}}=\mathbb{E}\sbr{e^{-\langle \lambda,W(t)\rangle}}\cdot\mathbb{E}\sbr{e^{-\langle \lambda,X(t)\rangle}},
$$
so the Laplace exponent $J_Z$ of $Z$ equals
\begin{equation} \label{LaplaceZ}
J_Z(\lambda)={{\frac{1}{2}\langle Q\lambda,\lambda\rangle+J_X(\lambda)}}.
\end{equation}

\begin{ex}[$\alpha$-stable martingales with $\bf{1<\alpha<2}$]\label{stable martingales}
A real valued stable martingale $Z^\alpha_t, t\geq 0$ with index $\alpha\in(1,2)$ and positive jumps only is a L\'evy process without Wiener part with L\'evy measure of the form
$$
\nu(\dd v):=\frac{1}{v^{\alpha+1}}\mathbf{1}_{\{v>0\}} \dd v.
$$ 
Its Laplace exponent is given by
\begin{align}\label{LAplace exp alfa stabilny}\nonumber
J_{Z^{\alpha}}(\lambda)&=\int_{0}^{+\infty}\left(e^{-\lambda v}-1+\lambda v\right)\frac{1}{v^{\alpha+1}} \dd v\\[1ex]
&=c_{\alpha} \lambda^\alpha, \quad \lambda \geq 0,
\end{align}
 with 
\begin{gather}\label{c alpha} 
c_\alpha:=\frac{\Gamma(2-\alpha)}{\alpha(\alpha-1)},
\end{gather}
where $\Gamma$ stands for the Gamma function.  Analogously one defines an $\alpha$-stable process with negative jumps only. 
\end{ex}
Note that the case of L\'evy martingale with the stability index $\alpha=2$ corresponds to the case when $Z^{\alpha}$ is a Wiener process without drift and with vanishing L\'evy measure.
\subsubsection{Projections of the noise}

For equation \eqref{rownanie 2} we consider the {\it projections} of $Z$ along $G$ given by
\begin{gather}\label{projection of Z}
Z^{G(x)}(t):=\langle G(x), Z(t)\rangle, \qquad x,t\geq 0.
\end{gather}
As linear transformations of $Z$, the projections form a family of L\'evy processes parametrized by $x\geq 0$. 
If $Z$ is a martingale, then $Z^{G(x)}$ is a real-valued L\'evy martingale for any $x\geq 0$.
It follows from the identity
$$
\mathbb{E}\sbr{e^{- \gamma\cdot Z^{G(x)}(t)}}=\mathbb{E}\sbr{e^{-\langle \gamma G(x), Z(t)\rangle}}, \quad \gamma\in\mathbb{R}, 
$$
and \eqref{LaplaceZ} that the Laplace exponent of $Z^{G(x)}$ equals 
\begin{gather}\label{Laplace ZG do uproszczenia}
J_{Z^{G(x)}}(\gamma)=J_Z(\gamma G(x))=\frac{1}{2}\gamma^2\langle Q G(x),G(x)\rangle+\int_{\mid y\mid>0}\left(e^{-\gamma \langle G(x),y\rangle}-1+\gamma\langle G(x),y\rangle\right)\nu(\dd y).
\end{gather}
Formula \eqref{Laplace ZG do uproszczenia} can be written in a simpler form by using the L\'evy measure $\nu_{G(x)}(\dd v)$ of $Z^{G(x)}$, which is the {\it image} 
of the L\'evy measure $\nu(dy)$ under the linear transformation $y\mapsto \langle G(x), y\rangle$. This measure is given by 
\begin{gather}\label{nu_G}
\nu_{G(x)}(A):=\nu \{y \in \R^d: \langle G(x),y\rangle\in A \} , \quad A\in\mathcal{B}(\mathbb{R}).
\end{gather}
From \eqref{Laplace ZG do uproszczenia}  we obtain that
\begin{gather}\label{Laplace ZG}
J_{Z^{G(x)}}(\gamma)=\frac{1}{2}\gamma^2\langle Q G(x),G(x)\rangle+\int_{\mid v\mid>0}\left(e^{-\gamma v}-1+\gamma v \right)\nu_{G(x)}(\dd v).
\end{gather}
Thus the characteristic triplet of the projection $Z^{G(x)}$ has the form
\begin{gather}\label{charakterystyki rzutu}
\left(-\int_{\mid v\mid>1}y \ \nu_{G(x)}(\dd v), \ \langle Q G(x),G(x)\rangle, \ \nu_{G(x)}(\dd v)\mid_{v\neq 0}\right).
\end{gather}
Above we used the restriction $\nu_{G(x)}(\dd v)\mid_{v\neq 0}$ by cutting off zero which may be an atom of $\nu_{G(x)}(\dd v)$.

\subsection{Preliminary characterization of generating equations}

In Proposition \ref{prop wstepny} below we provide a preliminary characterization for \eqref{rownanie 2} to be a generating equation. 
Note that the independence of coordinates of $Z$ is not assumed here. The central role here play the noise projections \eqref{projection of Z}. The result is deduced from Theorem 5.3 in \cite{FilipovicATS}, where the generator of a general non-negative Markovian  short rate process for affine models was characterized.

\begin{prop}\label{prop wstepny} 
Let $Z$ be a L\'evy martingale with characteristic triplet  \eqref{chrakterystyki Z} and $Z^{G(x)}$ be its projection \eqref{projection of Z} with the Le\'vy measure $\nu_{G(x)}(\dd v)$ given by \eqref{nu_G}.
\begin{enumerate}[(A)] 
\item Equation \eqref{rownanie 1} generates an affine model if and only if the following conditions are satisfied:
\begin{enumerate}[a)]
\item For each $x \ge 0$ the support of $\nu_{G(x)}$ is contained in $[0, +\ns)$ which means that $Z^{G(x)}$ has positive jumps only, i.e.  for each $t\geq 0$, with probability one,
\begin{gather}\label{Z^G positive jumps}
\triangle Z^{G(x)}(t):=Z^{G(x)}(t)-Z^{G(x)}(t-)=\langle G(x), \triangle Z(t)\rangle\geq 0.
\end{gather}
\item The jump part of $Z^{G(0)}$ has finite variation, i.e.
\begin{gather}\label{nu G0 finite variation}
\int_{(0,+\infty)}v \ \nu_{G(0)}(\dd v)<+\infty.
\end{gather}
\item The characteristic triplet \eqref{charakterystyki rzutu} of $Z^{G(x)}$ is linear in $x$, i.e.
\begin{align}\label{mult. CIR condition}
\frac{1}{2}\langle Q G(x), G(x)\rangle&=cx, \quad x\geq 0,\\[1ex]\label{rozklad nu G(x)}
\nu_{G(x)}(\dd v)\mid_{(0,+\infty)}&=\nu_{G(0)}(\dd v)\mid_{(0,+\infty)}+x\mu(\dd v), \quad x\geq 0,
\end{align}
for some $c\geq 0$ and a measure  $\mu(\dd v) \ \text{on} \ (0,+\infty) \ \text{satisfying}$
\begin{gather}\label{war calkowe na mu}
\int_{(0,+\infty)}(v \wedge v^2)\mu(\dd v)<+\infty.
\end{gather}
\item The function $F$ is affine, i.e.
\begin{gather}\label{linear drift}
F(x)=ax+b, \ \text{where} \ a\in\mathbb{R}, \ b\geq\int_{(1,+\infty)}(v-1)\nu_{G(0)}(\dd v) .
\end{gather}
\end{enumerate}
\item Equation \eqref{rownanie 1} generates an affine model if and only if the generator of $R$ is given by
\begin{align}\label{generator R w tw}\nonumber
\mathcal{A}f(x)=cx f^{\prime\prime}(x)&+\Big[ax +b+\int_{(1,+\infty)}(1 -v)\{\nu_{G(0)}(\dd v)+x\mu(\dd v)\}\Big]f^{\prime}(x)\\[1ex]
&+\int_{(0,+\infty)}[f(x+v)-f(x)-f^{\prime}(x)(1\wedge v)]\{\nu_{G(0)}(\dd v)+x\mu(\dd v)\}.
\end{align}
for $f\in\mathcal{L}(\Lambda)\cup C_c^2(\mathbb{R}_{+})$, where 
$\mathcal{L}(\Lambda)$ is the linear hull of $\Lambda:=\{f_\lambda:=e^{-\lambda x}, \lambda\in(0,+\infty)\}$
and $C_c^2(\mathbb{R}_{+})$ stands for the set of twice continuously differentiable functions with compact support in $[0,+\infty)$. The constants $a,b,c$ and the measures $\nu_{G(0)}(\dd v), \mu(\dd v)$ are those from part (A).
\end{enumerate}
\end{prop}

The poof of Proposition \ref{prop wstepny} is postponed to Appendix.

Note that conditions  \eqref{mult. CIR condition}-\eqref{rozklad nu G(x)} describe the distributions of the noise projections. In the sequel we use an equivalent formulation of  \eqref{mult. CIR condition}-\eqref{rozklad nu G(x)}  involving the Laplace exponents of \eqref{projection of Z}. Taking into account \eqref{Laplace ZG} we obtain the following.

\begin{rem}\label{rem warunki w eksp. Laplacea}
	The conditions \eqref{mult. CIR condition} and \eqref{rozklad nu G(x)}  are equivalent to the following decomposition of the Laplace exponent of $Z^G$:
	\begin{gather}\label{war na exp Laplaca}
	J_{Z^{G(x)}}(b)=cb^2x+J_{\nu_{G(0)}}(b)+x J_{\mu}(b), \quad b,x\geq 0,
	\end{gather}
	where
	\begin{gather}\label{def Jmu i Jnu0}
	J_{\mu}(b):=\int_{0}^{+\infty}(e^{-bv}-1+bv)\mu(\dd v), \quad J_{\nu_{G(0)}}(b):=\int_{0}^{+\infty}(e^{-bv}-1+bv)\nu_{G(0)}(\dd v).
	\end{gather}
\end{rem}

\subsubsection{Problem formulation}

In virtue of part $(A)$ of Proposition \ref{prop wstepny} we see that the drift $F$ of a generating equation is an affine function while 
the function $G$ and the noise $Z$ must provide projections $Z^{G(x)}, x\geq 0$ with particular distributions. 
Their characteristic triplets are characterized by a constant 
$c\geq 0$ carrying information on the variance of the Wiener part and two measures 
$\nu_{G(0)}(\dd v)$, $\mu(\dd v)$ describing jumps. 
A pair $(G,Z)$ for which the projections $Z^{G(x)}$ satisfy \eqref{nu G0 finite variation}-\eqref{war calkowe na mu}
will be called {\it a generating pair}. Note that the concrete forms of the measures $\nu_{G(0)}(\dd v)$, $\mu(\dd v)$
are, however, not specified. As for $Z$ with independent coordinates of infinite variation necessarily $G(0)=0$, see Proposition \ref{rem o G(0)=0}, and, consequently, $\nu_{G(0)}(\dd v)$ vanishes, our goal is to determine the measure $\mu(\dd v)$ in this case.

Having the required form of $\mu(\dd v)$ at hand one knows the distributions of the noise projections $Z^{G(x)}$ and, by  part $(B)$ of Proposition \ref{prop wstepny}, also the generator of the solution of \eqref{rownanie 2}. The generating pairs $(G,Z)$ can not be, however, 
uniquely determined, except the one-dimensional case.  This issue is discussed in 
Section \ref{Low-dimensional generating equations} and Section \ref{section Non-uniqueness in the multidimensional case} below.
For this reason we construct canonical representations - generating equations with noise projections corresponding to a given form of the measure $\mu(\dd v)$.

\subsubsection{One-dimensional generating equations}\label{Low-dimensional generating equations}

Let us summarize known facts on generating equations in the case $d=1$. If  $Z=W$ is a Wiener process, the only generating equation is the classical CIR equation
\begin{gather}\label{CIR equation}
\dd R(t)=(aR(t)+b)\dd t+C\sqrt{R(t)}\dd W(t), 
\end{gather}
with $a\in\mathbb{R}$, $b,C\geq 0$, see \cite{CIR}. The 
case with a general one-dimensional L\'evy process $Z$ was studied in \cite{BarskiZabczykCIR}, \cite{BarskiZabczyk} and
\cite{BarskiZabczykArxiv} with the following conclusion. If the variation of $Z$ is infinite and $G \not\equiv 0$, then $Z$ must be an $\alpha$-stable process with index $\alpha\in(1,2]$, with either positive or negative jumps only,  and \eqref{rownanie 1} has the form
\begin{gather}\label{CIR equation geenralized}
\dd R(t)=(aR(t)+b)\dd t+C\cdot R(t-)^{{1}/{\alpha}}\dd Z^{\alpha}(t),
\end{gather}
with $a\in\mathbb{R}, b\geq 0$ and $C$ such that it has the same sign as the jumps of $Z^\alpha$. Clearly, for $\alpha=2$ equation \eqref{CIR equation geenralized} becomes \eqref{CIR equation}. If $Z$ is of finite variation then the noise enters \eqref{rownanie 1} in the additive way, that is 
\begin{gather}\label{Vasicek equation geenralized}
\dd R(t)=(aR(t)+b)\dd t+C \ \dd Z(t).
\end{gather}
Here $Z$ can be chosen as an arbitrary process with positive jumps, $a\in\mathbb{R}, C\geq 0$ and 
$$
b\geq C \int_{0}^{+\infty}y \ \nu(\dd y),
$$
where $\nu(\dd y)$ stands for the L\'evy measure of $Z$.  The variation of $Z$ is finite, so is the right side above.
Recall, \eqref{Vasicek equation geenralized} with $Z$ being a Wiener process is the well known Vasi\v cek equation, see \cite{Vasicek}. Then the short rate is a Gaussian process, hence it takes negative values with positive probability. 
This drawback is eliminated by the jump version of  the Vasi\v cek equation \eqref{Vasicek equation geenralized}, where the solution never falls below zero.

It follows that the triplet $(c,\nu_{G(0)}(\dd v),\mu(\dd v))$ from  Proposition \ref{prop wstepny} takes for the equations above the following forms
\begin{enumerate} [a)]
\item $c \geq 0, \ \nu_{G(0)}(\dd v)\equiv 0, \ \mu(\dd v)\equiv 0$; \\[1ex]
This case corresponds to the classical CIR equation \eqref{CIR equation} where $c=\frac{1}{2}C^2$.
\item $c=0, \ \nu_{G(0)}(\dd v)\equiv0, \ \mu(\dd v)- \text{$\alpha$-stable}, \ \alpha\in(1,2)$;\\[1ex]
In this case \eqref{rownanie 2} becomes the generalized CIR equation with $\alpha$-stable noise \eqref{CIR equation geenralized}.
\item $c=0, \ \nu_{G(0)}(\dd v)- \text{any measure on $(0,+\infty)$ of finite variation}, \ \mu(\dd v)\equiv 0$;\\[1ex]
Here \eqref{rownanie 2} becomes the generalized Vasi\v cek equation \eqref{Vasicek equation geenralized}.
\end{enumerate}
Note the one to one correspondence between the triplets $(c,\nu_{G(0)}(\dd v),\mu(\dd v))$  and generating pairs $(G,Z)$ which holds up to multiplicative constants.

\subsubsection{Non-uniqueness in the multidimensional case}\label{section Non-uniqueness in the multidimensional case}

In the case $d>1$ one should not expect a one to one correspondence between the triplets $(c,\nu_{G(0)}(\dd v),\mu(\dd v))$ and the generating equations \eqref{rownanie 2}. The reason is that the distribution of the noise projections $Z^{G(x)}$
does not determine the pair $(G,Z)$ in a unique way. Our illustrating example below shows two different equations 
driven by L\'evy processes with independent coordinates which provide the same short rate $R$.

\begin{ex}\label{ex different eq ident sol}

Let us consider the following two equations 
\begin{align}\label{ex pierwsze rownanie}
\dd R(t)&=\langle G(R(t-)), \dd Z(t) \rangle, \quad R(0)=R_0,\quad t\geq 0,\\[1ex]\label{ex drugie rownanie}
d\bar{R}(t)&=\langle \bar{G}(\bar{R}(t-), \dd \bar{Z}(t))\rangle, \quad \bar{R}(0)=R_0,\quad t\geq 0,
\end{align}
where 
$$
G(x):=2^{-1/\alpha}\cdot (x^{1/ \alpha}, x^{1/ \alpha}), \quad Z:=(Z^\alpha_1, Z^\alpha_2), 
$$
and 
$$
\bar{G}(x):=(x^{1/ \alpha}, x^{1/ \alpha}), \quad \bar{Z}:=(\bar{Z}_1, \bar{Z}_2),
$$
with a fixed index $\alpha\in(1,2)$. We assume that the coordinates of $Z$ and $\bar{Z}$ are independent. Above $Z^{\alpha}_1, Z^{\alpha}_2$ stand for $\alpha$-stable martingales like in Example \ref{stable martingales} and $\bar{Z}_1, \bar{Z}_2$ are martingales with L\'evy measures
\[
\nu_1(\dd v) =  \frac{\dd v}{v^{\alpha +1}}\mathbf{1}_E(v), \quad \nu_2(\dd v) = \frac{\dd v}{v^{\alpha +1}} \mathbf{1}_{[0, +\ns) \setminus E}(v),
\]
respectively, where $E$ is a Borel subset of $[0, +\ns)$ such that 
$$
|E| = \int_{E} \dd v >0, \quad  \text{and} \quad |[0, +\ns) \setminus E| = \int_{[0, +\ns) \setminus E} \dd v  >0.
$$ 
The projections related to \eqref{ex pierwsze rownanie} and \eqref{ex drugie rownanie} take the forms
\begin{align*}
Z^{G(x)}(t)&=\langle G(x), Z(t)\rangle=x^{1/\alpha} 2^{-1/\alpha}(Z^\alpha_1(t)+Z^{\alpha}_2(t)), \quad x,t\geq 0,\\
\bar{Z}^{\bar{G}(x)}(t)&=\langle \bar{G}(x), \bar{Z}(t)\rangle=x^{1/\alpha} (\bar{Z}_1(t)+\bar{Z}_2(t)),\quad x,t\geq 0.
\end{align*}
Since both processes $2^{-1/\alpha} (Z^\alpha_1+Z^{\alpha}_2)$ and $\bar{Z}_1+\bar{Z}_2$ are 
$\alpha$-stable and have the same finite dimensional distributions, we obtain that 
$$
Z^{G(x)}=\bar{Z}^{\bar{G}(x)},
$$
in the sense of distribution. Moreover, the L\'evy measure of $Z^{G(x)}$ has the form
$$
x\cdot \frac{\dd v}{v^{\alpha+1}}\mathbf{1}_{\{v>0\}},  \quad x\geq 0,
$$
so it follows from \eqref{rozklad nu G(x)} that $(G,Z)$ is a generating pair
and that the solutions of \eqref{ex pierwsze rownanie} and \eqref{ex drugie rownanie} are identical.

Note that the triplet $(c,\nu_{G(0)},\mu(\dd v))$ from Proposition \ref{prop wstepny} is, for both pairs, of the form
$$
c=0, \ \nu_{G(0)}(\dd v)\equiv0, \ \mu(\dd v)- \text{$\alpha$-stable}, 
$$
so it coincides with the triplet $(b)$ in Section \ref{Low-dimensional generating equations}. Consequently, the solutions of 
\eqref{ex pierwsze rownanie} and \eqref{ex drugie rownanie} are the same as the solution of the equation
$$
dR(t)=(R(t-))^{1/\alpha} \dd Z^{\alpha}(t), \quad R(0)=R_0, \quad t\geq 0,
$$
with a one-dimensional $\alpha$-stable process $Z^\alpha$.

It follows, in particular, that the noise coordinates of a generating equation do not need to be stable processes.
\end{ex}

\section{Classification of generating equations}\label{section Classification of generating equations}
\subsection{Main results}\label{section Noise with independent coordinates}

This section deals with equation \eqref{rownanie 2} in the case when the coordinates of the martingale $Z$ are independent. In view of Proposition \ref{prop wstepny}  we are interested in characterizing possible distributions of projections $Z^G$ over all generating pairs $(G,Z)$.  By \eqref{Z^G positive jumps} the jumps of the projections are necessarily positive. As the coordinates of $Z$ are independent, they do not jump together. Consequently, we see that, for each $x\geq 0$ and $t \ge 0$
$$
\triangle Z^{G(x)}(t)=\langle G(x),\triangle Z(t)\rangle >0
$$
holds if and only if, for some $i=1,2,...,d$,
\begin{gather}\label{mucha}
G_i(x)\triangle Z_i(t)>0, \quad \triangle Z_j(t)=0, j\neq i.
\end{gather}
Condition \eqref{mucha} means that $G_i(x)$ and $\triangle Z_i(t)$ are of the same sign. We can consider only the case when both are positive, i.e.
$$
G_i(x)\geq 0, \quad i=1,2,...,d, \ x\geq 0, \qquad \triangle Z_i(t)\geq 0, \quad t> 0,
$$
because the opposite case can be turned into this one by replacing $(G_i,Z_i)$ with $(-G_i,-Z_i)$, $i=1,...,d$. The L\'evy measure $\nu_i(\dd y)$ of $Z_i$ is thus concentrated on $(0,+\infty)$ and, in view of \eqref{LaplaceZ}, the Laplace exponent of $Z_i$ takes the form
\begin{gather}\label{Laplace Zi}
J_i(b):=\frac{1}{2}q_{ii} b^2+\int_{0}^{+\infty}(e^{-b v}-1+b v)\nu_i(\dd v),\quad b \geq 0, \ i=1,2,...,d,
\end{gather}
with $q_{ii}\geq 0$. Recall, $q_{ii}$ stands on the diagonal of $Q$ - the covariance matrix of the Wiener part of $Z$.  
We will assume that $J_i, i=1,2,...,d$ are  {\it regularly varying at zero}. Recall, this means that
$$
\lim_{x\rightarrow 0^+}\frac{J_i(bx)}{J_i(x)}=\psi_i(b), \quad b> 0,\qquad i=1,2,...,d,
$$
for some function $\psi_i$. In fact $\psi_i$ is a power function, i.e.
$$
\psi_i(b)=b^{\alpha_i}, \quad b>0,
$$
with some $-\infty< \alpha_i<+\infty$ and $J_i$ is called to vary regularly with index $\alpha_i$. A characterization of regularly varying Laplace exponent in terms of the corresponding L\'evy measure is presented in Section \ref{sec slowly varying Laplace exp.}.

The distribution of noise projections are described by the following result.

\begin{tw} \label{TwNiez} Let $Z_1,...,Z_d$ be independent coordinates of the L\'evy martingale $Z$ in $\R^d$. Assume that $Z_1,...,Z_{d}$  satisfy
\begin{equation} \label{ass1}
\triangle Z_i(t)\geq 0 \text{ a.s.  for } t>0  \text{ and } Z_i \ \text{is of infinite variation}
\end{equation}
or 
\begin{equation} \label{ass2}
\triangle Z_i(t)\geq 0\text{ a.s.  for } t>0 \text{ and } G(0)=0.
\end{equation}
Further, let us assume that for all $i=1,\ldots, d$ the Laplace exponent  \eqref{Laplace Zi} of $Z_i$ varies regularly at zero and the components of the function  $G$ satisfiy
$$
G_i(x)\geq 0, \ x\in[0,+\infty), \quad G_i \ \text{is continuous on } [0,+\infty).
$$
Then \eqref{rownanie 2} generates an affine model if and only 
 if $F(x)=ax+b$, $a\in\mathbb{R}, b\geq 0$, and the Laplace exponent $J_{Z^{G(x)}}$ of $Z^{G(x)}=\langle G(x), Z\rangle$ is of the form 
\begin{gather}\label{postac J_ZG przy niezaleznych}
J_{Z^{G(x)}}(b) = x \sum_{k=1}^g\eta_{k}b^{\alpha_{{k}}}, \quad  \eta_{k}> 0, \quad \alpha_k\in(1,2],  \quad k=1,2,\ldots,g,
\end{gather}
with some $1 \le g \le d$ and $\alpha_k\neq\alpha_j$ for $k\neq j$.
\end{tw}

Theorem \ref{TwNiez} allows determining the form of the measure $\mu(\dd v)$ in Proposition \ref{prop wstepny}.

\begin{cor}\label{cor o postaci mu} 
Let the assumptions of Theorem \ref{TwNiez} be satisfied. If equation \eqref{rownanie 2} generates an affine model 
then the function $J_\mu$ defined in \eqref{def Jmu i Jnu0} takes the form
\begin{gather}\label{J mu postaaaccccc}
J_{\mu}(b)  = \sum_{k=l}^{g}\eta_{k}b^{\alpha_{{k}}}, \quad l\in \{1,2\}, \quad  \eta_{k}> 0, \quad \alpha_k\in(1,2),  \quad k=l,l+1,\ldots,g,
\end{gather}
with $1 \le g \le d$, $2>\alpha_l>...>\alpha_g>1$ (for the case $l=2, g=1$ we set $J_{\mu}\equiv 0$, which means that $\mu(\dd v)$ disappears). Above $l=2$ if $\alpha_1=2$ and $l=1$ otherwise. This means that $\mu(\dd v)$ is a weighted sum of $g+1-l$ stable measures with indices $\alpha_l,...,\alpha_g\in(1,2)$, i.e.
\begin{gather}\label{miary w klasach}
\mu(\dd v)=\tilde{\mu}(\dd v):= \frac{d_l}{v^{1+\alpha_l}}\mathbf{1}_{\{v>0\}}\dd v+...+\frac{d_g}{v^{1+\alpha_g}}\mathbf{1}_{\{v>0\}} \dd v,
 \end{gather}
with $d_i=\eta_i/ c_{\alpha_i},i=l,...,g$, where $c_{\alpha_i}$ is given by \eqref{c alpha} .
\end{cor}

Note that each generating equation can be identified by the numbers $a,b$ appearing in the formula for the function $F$ and $\alpha_1,...,\alpha_g; \eta_1,...,\eta_g$ from \eqref{postac J_ZG przy niezaleznych}. Since $\nu_{G(0)}(\dd v)=0$, see  Proposition \ref{rem o G(0)=0} in the sequel, the related generator of $R$ takes, by \eqref{generator R w tw}, the form
\begin{align}\label{generator niezalezne}\nonumber
 \mathcal{A}f(x)=cx f^{\prime\prime}(x)&+\Big[x\Big(a+\int_{(1,+\infty)}(1 -v)x\tilde{\mu}(\dd v)\Big)+b\Big]f^{\prime}(x)\\[1ex]
 &+\int_{(0,+\infty)}[f(x+v)-f(x)-f^{\prime}(x)(1\wedge v)]x\tilde{\mu}(\dd v),
\end{align}
with $\tilde{\mu}$ in \eqref{miary w klasach}. Recall, the constant $c$ above comes from the condition 
\begin{gather}\label{FG w klasach}
\frac{1}{2}\langle Q G(x),G(x)\rangle=cx, \quad  \quad x\geq 0,
\end{gather}
and, in view of Remark \ref{rem warunki w eksp. Laplacea}, $c=\eta_1$ if $\alpha_1=2$ and $c=0$ otherwise.
The class of processes with generator of the form \eqref{generator niezalezne} will be denoted by
\begin{gather}\label{klasa Ag}
\mathbb{A}_g(a,b;\alpha_1,\alpha_2,...,\alpha_g; \eta_1,...,\eta_g),
\end{gather}
 All generating equations with $d$-dimensional noise $Z$ satisfying assumptions of Theorem \ref{TwNiez} are thus splitted into $d$ disjoint subfamilies providing different short rates. Any two equations from \eqref{klasa Ag} with fixed parameters provide the same short rate, hence the same bond prices. For any class \eqref{klasa Ag} we construct below a {\it canonical representation}, which is an equation with the generator required in \eqref{klasa Ag} but with reduced noise dimension from $d$ to $g$ and stable noise coordinates. This construction allows interpreting the parameter $g$ in \eqref{klasa Ag} as a minimal number of random factors necessary to obtain the short rate corresponding to \eqref{klasa Ag} and $\alpha_1,\alpha_2,...,\alpha_g$ are the stability indices of the noise coordinates. This idea of classifying is similar to that of Dai and Singleton applied for multi-factor affine short rates in \cite{DaiSingleton}.

\begin{prop}[Canonical representation of $\mathbb{A}_g(a,b;\alpha_1,\alpha_2,...,\alpha_g; \eta_1,...,\eta_g)$]\label{prop canonical representation}
Let $R$ be the solution of \eqref{rownanie 2} with $F,G,Z$ satisfying the assumptions of Theorem \ref{TwNiez}. 
 Let $\tilde{Z}=(\tilde{Z}^{\alpha_1}_1,\tilde{Z}^{\alpha_2}_2,...,\tilde{Z}^{\alpha_g}_g)$ be a L\'evy martingale with independent stable coordinates
with indices $\alpha_k, k=1,2,...,g$, respectively, and $\tilde{G}(x)=(d^{1/\alpha_1}_1 x^{1/\alpha_1},...,d^{1/\alpha_g}_g x^{1/\alpha_g})$, $x \ge 0$, where $d_k:=\eta_k/c_{\alpha_k}$ and $c_{\alpha_k}$  are given by \eqref{c alpha}, $k=1,2,...,g$.
Then 
$$
J_{Z^{G(x)}}(b)=J_{\tilde{Z}^{\tilde{G}(x)}}(b), \quad b,x\geq 0.
$$
Consequently, if $\tilde{R}$ is the solution of the equation 
\begin{gather}\label{rownanie sklajane}
\dd \tilde{R}(t)=(a\tilde{R}(t)+b) \dd t+\sum_{k=1}^{g}d_k^{1/{\alpha_k}} \tilde{R}(t-)^{1/{\alpha_k}}\dd \tilde{Z}_k(t),
\end{gather}
then the generators of $R$ and $\tilde{R}$ are equal.
\end{prop}

Equation \eqref{rownanie sklajane} will be called the {\it canonical representation} of the class $\mathbb{A}_g(a,b;\alpha_1,\alpha_2,...,\alpha_g; \eta_1,...,\eta_g)$.
\vskip2ex
\noindent
{\bf Proof:} By \eqref{postac J_ZG przy niezaleznych} we need to show that 
\begin{gather*}\label{rozklad rzutu konstrukcja}
J_{\tilde{Z}^{\tilde{G}(x)}}(b)=x\sum_{k=1}^{g}\eta_k b^{\alpha_k}, \quad b,x\geq 0.
\end{gather*}
Recall, the Laplace exponent of $\tilde{Z}^{\alpha_k}_k$ equals
$J_k(b)=c_{\alpha_k}b^{\alpha_k}, k=1,2,...,g$. By independence and the form of $\tilde{G}$ we have
\begin{align*}
J_{\tilde{Z}^{\tilde{G}(x)}}(b)&=\sum_{k=1}^{g}J_k(b\tilde{G}_k(x))=\sum_{k=1}^{g}c_{\alpha_k} b^{\alpha_k}d_kx=x\sum_{k=1}^{g}\eta_k b^{\alpha_k}, \quad b,x\geq 0,
\end{align*}
as required. The second part of the thesis follows from Proposition \ref{prop wstepny}(B).\hfill$\square$

\vskip1ex
Clearly, in the case $d=1$ the noise dimension can not be reduced, so $g=d=1$ and $\mathbb{A}_1(a,b; 2;\eta_1)$ corresponds to the classical CIR equation \eqref{CIR equation} while $\mathbb{A}_1(a,b; \alpha;\eta_1), \alpha\in(1,2)$ to its generalized version 
\eqref{CIR equation geenralized}. Both classes are singletons and \eqref{CIR equation}, \eqref{CIR equation geenralized} are their canonical representations. The alpha-CIR equation from \cite{JiaoMaScotti} is a canonical representation of the class $\mathbb{A}_2(a,b; 2,\alpha; \eta_1,\eta_2)$ with $\alpha\in(1,2)$.

\subsubsection{Proofs}
The proofs of Theorem \ref{TwNiez} and Corollary \ref{cor o postaci mu} are preceded by two auxiliary results, i.e. Proposition \ref{bounds_alpha} and
Proposition \ref{rem o G(0)=0}. The first one provides some useful estimation for the function 
\begin{gather}\label{J}
J_{\rho}(b):=\int_{0}^{+\infty}(e^{-bv}-1+bv)\rho(\dd v), \quad b\geq 0,
\end{gather}
where the measure $\rho(\dd v)$ on $(0,+\ns)$ satisfies
\begin{gather}\label{nuJ}
0 < \int_0^{+\ns} \rbr{v^2\wedge v} \rho\rbr{\dd v} < +\ns.
\end{gather}
The second result shows that if all components of $Z$ are of infinite variation then $G(0)=0$.

\begin{prop} \label{bounds_alpha}
Let $J_{\rho}$ be a function given by \eqref{J} where the measure  $\rho$ satisfies \eqref{nuJ}. Then the function
$
(0,+\ns) \ni b \mapsto {J_{\rho}(b)}/{b}$ is strictly increasing and $\lim_{b \ra 0+}J_{\rho}(b)/b = 0$, while the function $(0,+\ns) \ni b \mapsto {J_{\rho}(b)}/{b^2}$
is strictly decreasing and $\lim_{b \ra +\ns}J_{\rho}(b)/b^2 = 0$. This yields, in particular, that, for any $b_0 >0$, 
\begin{gather}\label{oszacowania dwustronne J}
\frac{J_{\rho}\rbr{b_0}}{b_0^2}b^2 < J_{\rho}(b) < \frac{J_{\rho}\rbr{b_0}}{b_0}b, \quad b\in \rbr{0, b_0}.
\end{gather}
\end{prop} 
\noindent
{\bf Proof:} Let us start from the observation that the function 
$$
t \mapsto \frac{(1-e^{-t})t}{e^{-t}-1+t}, \quad t\geq 0,
$$
is strictly decreasing, with limit $2$ at zero and $1$ at infinity. This implies
 \begin{equation} \label{oszH}
 (e^{-t}-1+t) < (1-e^{-t})t < 2 (e^{-t}-1+t), \quad t \in (0, +\ns),
 \end{equation}
and, consequently,
$$
\int_{0}^{+\infty}(e^{-bv}-1+bv)\rho(\dd v) < \int_{0}^{+\infty}(1-e^{-bv})bv\ \rho(\dd v) < 2\int_{0}^{+\infty}(e^{-bv}-1+bv)\rho(\dd v), \quad b >0.
$$
This means, however, that
$$
J_{\rho}(b) < bJ_{\rho}^\prime(b) < 2J_{\rho}(b), \quad b > 0.
$$
So, we have
$$
\frac{1}{b} < \frac{J_{\rho}^\prime(b)}{J_{\rho}(b)}=\frac{d}{db}\ln J_{\rho}(b) < \frac{2}{b}, \quad b>0,
$$
and integration over some interval $[b_1,b_2]$, where $b_2 > b_1>0$, yields
$$
\ln b_2 - \ln b_1  < \ln J_{\rho}\rbr{b_2}-\ln J_{\rho}\rbr{b_1} < 2 \ln b_2 - 2 \ln b_1
$$
which gives that 
$$
\frac{J_{\rho}\rbr{b_2}}{b_2} > \frac{J_{\rho}\rbr{b_1}}{b_1}, \quad \frac{J_{\rho}\rbr{b_2}}{b_2^2} < \frac{J_{\rho}\rbr{b_1}}{b_1^2}.
$$

To see that $\lim_{b \ra 0+} {J_{\rho}\rbr{b}}/{b} = 0$ it is sufficient to use de l'H\^opital's rule,  \eqref{nuJ} and dominated convergence
$$
\lim_{b \ra 0+} \frac{J_{\rho}\rbr{b}}{b} = \lim_{b \ra 0+} {J'_{\rho}\rbr{b}} = \lim_{b \ra 0+}  \int_{0}^{+\infty}(1-e^{-bv}) v\ \rho(\dd v) = 0.
$$

To see that $\lim_{b \ra +\ns} {J_{\rho}\rbr{b}}/{b^2} = 0$ we also use de l'H\^opital's rule,  \eqref{nuJ} and dominated convergence.
If $\int_0^{+\ns} v\ \rho\rbr{\dd v} < +\ns$, then we have
$$
\lim_{b \ra +\ns} \frac{J_{\rho}\rbr{b}}{b^2} = \lim_{b \ra +\ns} \frac{J_{\rho}'\rbr{b}}{2b} = \frac{\int_0^{+\ns} v \rho\rbr{\dd v}}{+\ns}= 0.
$$
If $\int_0^{+\ns} v\ \rho\rbr{\dd v} = +\ns$ then we apply de l'H\^opital's rule twice and obtain
$$
\lim_{b \ra +\ns} \frac{J_{\rho}\rbr{b}}{b^2} = \lim_{b \ra +\ns} \frac{J_{\rho}'\rbr{b}}{2b} = \lim_{b \ra +\ns} \frac{J_{\rho}''\rbr{b}}{2} = \frac{1}{2} \lim_{b \ra +\ns} \int_{0}^{+\infty}e^{-bv} v^2\ \rho(\dd v) = 0.
$$
\hfill $\square$

\begin{prop}\label{rem o G(0)=0}
If $(G,Z)$ is a generating pair and all components of $Z$ are of infinite variation then $G(0)=0$. 
\end{prop}
\noindent 
{\bf Proof:}  Let $(G,Z)$ be a generating pair. Since the components of $Z$ are independent, its characteristic triplet \eqref{chrakterystyki Z} is such that $Q=\{q_{i,j}\}$ is a diagonal matrix, i.e.
$$
q_{ii}\geq 0, \quad q_{i,j}=0, \qquad i\neq j, \quad i,j=1,2,...,d,
$$
and the support of $\nu(\dd y)$ is contained in the positive half-axes of $\mathbb{R}^d$, see \cite{Sato} p.67. 
On the $i^{th}$ positive  half-axis 
\begin{gather}\label{rozlozenie ny}
\nu(\dd y)=\nu_i(dy_i),\qquad y=(y_1,y_2,...,y_d),
\end{gather}
for $i=1,2,...,d$.
 The $i^{th}$ coordinate of $Z$ is of infinite variation if and only if its Laplace exponent \eqref{Laplace Zi} is such that $q_{ii}>0$ or
\begin{gather}\label{rerere}
\int_{0}^{1}y_i\nu_i(\dd y_i)=+\infty,
\end{gather}
see \cite[Lemma 2.12]{Kyprianou}. It follows from \eqref{mult. CIR condition} that
$$
\frac{1}{2}\langle QG(x), G(x)\rangle=\frac{1}{2}\sum_{j=1}^{d}q_{jj}G_j^2(x)=cx,
$$
so if $q_{ii}>0$ then $G_i(0)=0$. If it is not the case, using \eqref{rozlozenie ny} and  \eqref{nu G0 finite variation} we see that the integral
\begin{align*}
\int_{(0,+\infty)}v\nu_{G(0)}(\dd v) & =\int_{\mathbb{R}^d_{+}}\langle G(0),y\rangle \nu(\dd y) \\ 
&  =\sum_{j=1}^{d}\int_{(0,+\infty)}G_j(0)y_j \ \nu_j(\dd y_j) =\sum_{j=1}^{d}
G_j(0) \int_{(0,+\infty)}y_j \ \nu_j(\dd y_j),
\end{align*}
is finite, so if  \eqref{rerere} holds then $G_i(0)=0$.\hfill $\square$

 \vskip2ex

\noindent
{\bf Proof of Theorem \ref{TwNiez}:} By assumption \eqref{ass1} and Proposition \ref{rem o G(0)=0} or by assumption \eqref{ass2}  we have $G(0)=0$, so it follows from Remark \ref{rem warunki w eksp. Laplacea} that
\begin{equation} \label{eq:dwa_zero}
J_{Z^{G(x)}}(b) = J_{1}(bG_{1}(x))+J_2(bG_2(x))+...+J_d(bG_d(x))=x\tilde{J}_{\mu}(b), \quad b,x\geq 0,
\end{equation} 
where $\tilde{J}_{\mu}(b) = c b^2 + {J}_{\mu}(b)$, $c\ge 0$ and ${J}_{\mu}(b)$ is given by \eqref{def Jmu i Jnu0}.
This yields
\begin{equation}
\frac{J_{1}\rbr{b\cdot G_{1}(x)}}{J_{1}\rbr{G_{1}(x)}}\cdot\frac{J_{1}\rbr{G_{1}(x)}}{x}+\ldots+\frac{J_{d}\rbr{b\cdot G_{d}(x)}}{J_{d}\rbr{G_{d}(x)}}\cdot\frac{J_{d}\rbr{G_{d}(x)}}{x}=\tilde{J}_{\mu}(b),\label{eq:dwa}
\end{equation}
where in the case $G_i(x) = 0$ we set $\frac{J_{i}\rbr{b\cdot G_{i}(x)}}{J_{i}\rbr{G_{i}(x)}}\cdot\frac{J_{i}\rbr{ G_{i}(x)}}{x} = 0$. Without loss of generality we may assume that $J_{1}$, $J_{2}$,$\ldots$,$J_{d}$ are non-zero (thus positive for positive arguments).  
By assumption, $J_{i}$, $i=1,2,\ldots, d$ vary regularly at $0$ with some indices $\alpha_{i}$, $i=1,2,\ldots,d$, 
so for $b>0$
\begin{equation} \label{eq:trzy}
\lim_{y\ra0+}\frac{J_{i}\rbr{b\cdot y}}{J_{i}(y)}=b^{\alpha_{i}}.
\end{equation}
Assume that 
\[
\alpha_{1}=\ldots=\alpha_{i\rbr{1}}>\alpha_{i\rbr{1}+1}=\ldots=\alpha_{i\rbr{2}}>\ldots\ldots>\alpha_{i\rbr{g-1}+1}=\ldots=\alpha_{i\rbr{g}}=\alpha_{d},
\]
where $i(g) = d$. Let us denote $i_0 = 0$ and
\begin{equation} \label{limits}
\eta_{k}(x) :=\frac{J_{i\rbr{k-1}+1}\rbr{G_{i\rbr{k-1}+1}(x)}+\ldots+J_{i\rbr{k}}\rbr{G_{i\rbr{k}}(x)}}{x}, \quad k=1,2,\ldots, g.
\end{equation}
We can rewrite equation \eqref{eq:dwa} in the form
\begin{equation}
\sum_{k=1}^{g} \rbr{\sum_{i=i\rbr{k-1}+1}^{i\rbr{k}} \frac{J_{i}\rbr{b\cdot G_{i}(x)}}{J_{i}\rbr{G_{i}(x)}}\cdot\frac{J_{i}\rbr{ G_{i}(x)}}{x}}=\tilde{J}_{\mu}(b).\label{eq:trzyy}
\end{equation}
By passing to the limit as $x\ra0+$, from \eqref{eq:trzy} and \eqref{eq:trzyy} we get 
\begin{align}
b^{\alpha_{i\rbr{1}}}\rbr{ \lim_{x \ra 0+} \eta_{1}(x)} +\ldots+   b^{\alpha_{i\rbr{g}}} \rbr{\lim_{x \ra 0+} \eta_{g}(x)}  =\tilde{J}_{\mu}(b), \label{eq:trzyyy}
\end{align} 
thus
\begin{gather}\label{J mu tilde sum power}
\tilde{J}_{\mu}(b)=\sum_{k=1}^g \eta_{k}b^{\alpha_{i\rbr{k}}},
\end{gather}
provided that the limits $\eta_{k} := \lim_{x \ra 0+} \eta_{k}(x)$, $k=1,2,\ldots, g$, exist. 
Thus it remains to prove that for  $k=1,2,\ldots, g$ the limits $\lim_{x \ra 0+} \eta_{k}(x)$ indeed exist {and that $\alpha_{i(k)} \in (1,2]$.}

First we will prove that $\lim_{x \ra 0+} \eta_{g}(x)$ exists.
Assume, by contrary, that this is not true, so 
\begin{equation} 
\limsup_{x \ra 0+} \eta_{g}(x) - \liminf_{x \ra 0+} \eta_{g}(x) \ge \delta >0.
\label{sequencess}
\end{equation}
It follows from \eqref{eq:dwa_zero} that
\begin{equation} \label{ogr}
\frac{J_1(G_1(x))+J_2(G_2(x))+...+J_d(G_d(x))}{x} = \sum_{k=1}^g \eta_k(x)=\tilde{J}_{\mu}(1).
\end{equation}
Let now $b_0 \in (0,1)$ be small enough so that 
\begin{equation} \label{oszacowanie}
\tilde{J}_{\mu}(1) b_0^{\alpha_{i\rbr{g-1}} - \alpha_{i(g)}} < \frac{\delta}{6}.
\end{equation}
Let us set in \eqref{eq:trzyy} $b=b_0$ and then divide both sides of \eqref{eq:trzyy} by $b_0^{\alpha_{i(g)}}$. 
It follows from \eqref{ogr} that each term $\frac{J_{i}\rbr{ G_{i}(x)}}{x}$, $i=1,2,\ldots,d$, is bounded by 
 $\tilde{J}_{\mu}(1)$. From this and \eqref{eq:trzy} for  $x>0$ sufficiently close to $0$ we have 
\[
\eta_g(x) - \frac{\delta}{6} \le \frac{1}{b_0^{\alpha_{i(g)}}} \rbr{ \sum_{i=i\rbr{g-1}+1}^{i\rbr{g}} \frac{{J}_{i}\rbr{b_0\cdot G_{i}(x)}}{J_{i}\rbr{G_{i}(x)}}\cdot\frac{J_{i}\rbr{ G_{i}(x)}}{x}} \le \eta_g(x) + \frac{\delta}{6}
\]
and 
\begin{align*}
\frac{1}{b_0^{\alpha_{i(g)}}}  \sum_{k=1}^{g-1}  \rbr{ \sum_{i=i\rbr{k-1}+1}^{i\rbr{k}} \frac{J_{i}\rbr{b_0\cdot G_{i}(x)}}{J_{i}\rbr{G_{i}(x)}}\cdot\frac{J_{i}\rbr{ G_{i}(x)}}{x}}  \le \sum_{k=1}^{g-1}  2 b_0^{\alpha_{i(k)} - \alpha_{i(g)}} \eta_k(x) \\
\le 2 b_0^{\alpha_{i(g-1)} - \alpha_{i(g)}} \tilde{J}_{\mu}(1)
\end{align*}
thus from \eqref{eq:trzyy}, two last estimates and \eqref{oszacowanie}
\[
\eta_g(x) - \frac{\delta}{6} \le \frac{\tilde{J}_{\mu}(b_0)}{b_0^{\alpha_{i(g)}}} \le \eta_g(x) + \frac{\delta}{6}+ 2\tilde{J}_{\mu}(1) b_0^{\alpha_{i(g-1)} - \alpha_{i(g)}} < \eta_g(x) + \frac{\delta}{2}.
\]
But this contradicts \eqref{sequencess}
since we must have 
\[
 \limsup_{x \ra 0+} \eta_g(x) \le \frac{\tilde{J}_{\mu}(b_0)}{b_0^{\alpha_{i(g)}}} + \frac{\delta}{6}, \quad \liminf_{x \ra 0+} \eta_g(x) \ge \frac{\tilde{J}_{\mu}(b_0)}{b_0^{\alpha_{i(g)}}} - \frac{\delta}{2}.
\]

Having proved the existence of the limits $\lim_{x \ra 0+} \eta_{g}(x)$, ..., $\lim_{x \ra 0+} \eta_{g-m+1}(x)$ we can proceed similarly to prove the existence of the limit $\lim_{x \ra 0+} \eta_{g-m}(x)$. 
Assume that $\lim_{x \ra 0+} \eta_{g-m}(x)$ does not exist, so \begin{equation} 
\limsup_{x \ra 0+} \eta_{g-m}(x) - \liminf_{x \ra 0+} \eta_{g-m}(x) \ge \delta >0.
\label{sequencess1}
\end{equation}
Let $b_0 \in (0,1)$ be small enough so that 
\begin{equation} \label{oszacowanie1}
\tilde{J}_{\mu}(1) b_0^{\alpha_{i\rbr{g-m-1}} - \alpha_{i(g-m)}} < \frac{\delta}{8}.
\end{equation}
Let us set in \eqref{eq:trzyy} $b=b_0$ and then divide both sides of \eqref{eq:trzyy} by $b_0^{\alpha_{i(g-m)}}$. 
For  $x>0$ sufficiently close to $0$ we have 
\[
\eta_{g-m}(x) - \frac{\delta}{8} \le \frac{1}{b_0^{\alpha_{i(g-m)}}} \sum_{i=i\rbr{g-m-1}+1}^{i\rbr{g-m}} \frac{{J}_{i}\rbr{b_0\cdot G_{i}(x)}}{J_{i}\rbr{G_{i}(x)}}\cdot\frac{J_{i}\rbr{ G_{i}(x)}}{x} \le \eta_{g-m}(x) + \frac{\delta}{8},
\]
\begin{align*}
\frac{1}{b_0^{\alpha_{i(g-m)}}}  \sum_{k=1}^{g-m-1}  \rbr{ \sum_{i=i\rbr{k-1}+1}^{i\rbr{k}} \frac{J_{i}\rbr{b_0\cdot G_{i}(x)}}{J_{i}\rbr{G_{i}(x)}}\cdot\frac{J_{i}\rbr{ G_{i}(x)}}{x}}  \le \sum_{k=1}^{g-m-1}  2 b_0^{\alpha_{i(k)} - \alpha_{i(g-m)}} \eta_k(x) \\
\le 2 b_0^{\alpha_{i(g-m-1)} - \alpha_{i(g-m)}} \tilde{J}_{\mu}(1)
\end{align*}
and
\begin{align*}
\sum_{k=g-m+1}^{g}  \frac{b_0^{\alpha_{i(k)}} \eta_k}{b_0^{\alpha_{i(g-m)}}} - \frac{\delta}{8} & \le \frac{1}{b_0^{\alpha_{i(g-m)}}} \sum_{k=g-m+1}^{g}  \sum_{i=i\rbr{k-1}+1}^{i\rbr{k}} \frac{{J}_{i}\rbr{b_0\cdot G_{i}(x)}}{J_{i}\rbr{G_{i}(x)}}\cdot\frac{J_{i}\rbr{ G_{i}(x)}}{x} \\
& \le \sum_{k=g-m+1}^{g}  \frac{b_0^{\alpha_{i(k)}} \eta_k}{b_0^{\alpha_{i(g-m)}}}  + \frac{\delta}{8}
\end{align*}
thus from \eqref{eq:trzyy}, last three estimates and \eqref{oszacowanie1}
\begin{align*}
\eta_{g-m}(x) - \frac{\delta}{4} & \le \frac{J_{\mu}(b_0)}{b_0^{\alpha_{i(g-m)}}} - \sum_{k=g-m+1}^{g}  \frac{b_0^{\alpha_{i(k)}} \eta_k}{b_0^{\alpha_{i(g-m)}}} \\
 & \le \eta_{g-m}(x) + \frac{\delta}{4}+ 2\tilde{J}_{\mu}(1) b_0^{\alpha_{i(g-1)} - \alpha_{i(g)}} <\eta_{g-m}(x) +  \frac{\delta}{2}.
\end{align*}
But this contradicts \eqref{sequencess1}.

Now we are left with the proof that for $k=1,2,\ldots,g$, $\alpha_{i(k)} \in (1,2]$. Since the Laplace exponent of $Z_i$ is given by 
\eqref{Laplace Zi}, by Proposition \ref{bounds_alpha} we necessarily have that $J_i$ varies regularly with index $\alpha_i \in [1,2], i=1,2,...,d$. Thus it remains to prove that $\alpha_i >1, i=1,2,...,d$. If it was not true we would have $\alpha_{i(g)}=1$ in \eqref{J mu tilde sum power} and $\eta_g>0$. Then
$$
\lim_{b\ra 0+} \tilde{J}_{\mu}(b)/b = \lim_{b\ra 0+} J_{\mu}(b)/b =\eta_{g}>0,
$$ 
but, again, by Proposition \ref{bounds_alpha} it is not possible. 
\hfill $\square$

\vskip2ex
\noindent
{\bf Proof of Corollary \ref{cor o postaci mu} :} From Remark \ref{rem warunki w eksp. Laplacea} and Theorem \ref{TwNiez} we know that 
\[
J_{Z^{G(x)}}(b) =  x  c b^2 + x{J}_{\mu}(b) = x \sum_{k=1}^g\eta_{k}b^{\alpha_{{k}}}, 
\]
where $1\le g \le d$, $\eta_{k}> 0$, $\alpha_k\in(1,2]$, $\alpha_k\neq\alpha_j$, $k,j=1,2,\ldots,g$, $c\ge 0$.
Without loss of generality we may assume that $2 \ge \alpha_1 > \alpha_2 > \ldots > \alpha_g >1$. Thus, since the Laplace exponent is nonnegative, $ x{J}_{\mu}(b)$ is of the form
\begin{gather}\label{pierwsza postac}
x{J}_{\mu}(b) = x\sum_{k=1}^g\eta_{k}b^{\alpha_{{k}}}, \qquad \text {if} \ c=0,
\end{gather}
or
\begin{gather}\label{druga postac}
x{J}_{\mu}(b) = x\sbr{(\eta_1-c)b^2+ \sum_{k=2}^g\eta_{k}b^{\alpha_{{k}}}}, \qquad \text {if} \ 0<c\leq\eta_1 \ \text {and} \ \alpha_1=2.
\end{gather}
In the case \eqref{pierwsza postac} we need to show that $\alpha_1<2$. If it was not true, we would have
\[
\lim_{b \ra +\ns} \frac{{J}_{\mu}(b)}{b^2} = \eta_1 >0,
\]
but this contradicts Proposition \ref{bounds_alpha}. In the same way we prove that $\eta_1=c$ in \eqref{druga postac}.
This proves the required representation \eqref{J mu postaaaccccc}.
\hfill $\square$

\subsection{Characterization of regularly varying Laplace exponents}\label{sec slowly varying Laplace exp.}

In this section we reformulate the assumption that $J_i, i=1,...,d$, vary regularly at zero 
in terms of the behaviour of the L\'evy measures of $Z_i, i=1,...,d$. As our considerations are componentwise,
we write for simplicity $\nu(\dd v):=\nu_i(\dd v)$ for the L\'evy measure of $Z_i$ and $J:=J_i$ for its Laplace exponent.

\begin{prop}\label{prop J ragularly varying warunki na nu}

Let $\nu(\dd v)$ be such that 
\begin{gather}\label{war calk w prop charakt.}
\int_{0}^{+\infty}(y^2\wedge y) \ \nu(dy)<+\infty.
\end{gather}
Let $\tilde{\nu}(\dd v)$ be the measure
$$
\tilde{\nu}(\dd v):=v^2\nu(\dd v),
$$
and $\tilde{F}$ its cumulative distribution function, i.e.
$$
\tilde{F}(v):=\tilde{\nu}((0,v))=\int_{0}^{v}u^2\nu(\dd u),\quad v\geq 0. 
$$
Then, for $\alpha\in(1,2)$, the following conditions are equivalent
\begin{gather}\label{J varying regggg.}
\lim_{x\rightarrow 0^{+}}\frac{J(bx)}{J(x)}=b^{\alpha}, \quad b\geq 0,
\end{gather}
$$
\lim_{y\rightarrow +\infty}\frac{\tilde{F}(by)}{\tilde{F}(y)}=b^{2-\alpha}, \quad b\geq 0.
$$
If, additionally, $\nu(\dd v)$ has a density function $g(v)$ such that 
\begin{gather}\label{niecalkowalnosc y2}
\int_{0}^{+\infty}v^2 g(v) \nu(\dd v)=+\infty,
\end{gather}
then \eqref{J varying regggg.} is equivalent 
to the condition
$$
\lim_{y\rightarrow +\infty}\frac{g(by)}{g(y)}=b^{-\alpha-1}, \quad b > 0.
$$
\end{prop}
{\bf Proof:}  Under \eqref{war calk w prop charakt.} the function $J$ given by \eqref{J} is well defined for $b\geq 0$, twice differentiable and 
\begin{gather*}
J^\prime(b)=\int_{0}^{+\infty}v(1-e^{-bv})\nu(\dd v), \quad 
J^{\prime\prime}(b)=\int_{0}^{+\infty}v^2e^{-bv}\nu(\dd v), \quad b\geq 0,
\end{gather*}
see \cite{Rusinek}, Lemma 8.1 and Lemma 8.2.
This implies that
\begin{align*}
\lim_{x\rightarrow 0^{+}}\frac{J(bx)}{J(x)}&=b\cdot \lim_{x\rightarrow 0^{+}}\frac{J^\prime(bx)}{J^\prime(x)}=
b^2\cdot \lim_{x\rightarrow 0^{+}}\frac{J^{\prime\prime}(bx)}{J^{\prime\prime}(x)}\\[1ex]
&=b^2\cdot \lim_{x\rightarrow 0^{+}}\frac{\int_{0}^{+\infty}e^{-bxv}v^2\nu(\dd v)}{\int_{0}^{+\infty}e^{-xv}v^2\nu(\dd v)}.
\end{align*}
Consequently, by \eqref{J varying regggg.} 
\begin{gather}\label{do TAubera}
\lim_{x\rightarrow 0^{+}}\frac{\int_{0}^{+\infty}e^{-bxv}v^2\nu(\dd v)}{\int_{0}^{+\infty}e^{-xv}v^2\nu(\dd v)}=b^{\alpha-2}.
\end{gather}
Notice, that the left side  is a quotient of two transforms of the measure $\tilde{\nu}(\dd v)$.
By the Tauberian theorem, see Theorem 1, Sec. XIII.5 in \cite{Feller}, we have that \eqref{do TAubera} holds if and only if
$$
\frac{\tilde{F}(by)}{\tilde{F}(y)}\underset{y\rightarrow +\infty}{\longrightarrow}b^{2-\alpha}, \quad b\geq 0.
$$
If $\nu(\dd v)$ has a density $g(v)$ satisfying \eqref{niecalkowalnosc y2} then
\begin{align*}
\lim_{y\rightarrow +\infty}\frac{\tilde{F}(by)}{\tilde{F}(y)}&=\lim_{y\rightarrow +\infty}\frac{\int_{0}^{by}u^2g(u)\dd u}{\int_{0}^{y}u^2g(u)\dd u}
=\lim_{y\rightarrow +\infty}\frac{b\cdot (by)^2g(by)}{y^2g(y)}\\[1ex]
&=b^3\cdot\lim_{y\rightarrow +\infty}\frac{g(by)}{g(y)}.
\end{align*}
It follows that
$$
\lim_{y\rightarrow +\infty}\frac{g(by)}{g(y)}=b^{-\alpha-1}.
$$
which proves the result.\hfill$\square$

\begin{rem}
By general characterization of regularly varying functions we see that the functions
$\tilde{F}$ and $g$ from Proposition \ref{prop J ragularly varying warunki na nu} must be of the forms
$$
\tilde{F}(b)=b^{2-\alpha}L(b), \quad b\geq 0,
$$
$$
g(b)=b^{-\alpha-1}\tilde{L}(b), \quad b\geq 0,
$$
where $L$ and $\tilde{L}$ are slowly varying functions at $+\infty$, i.e.
$$
\frac{L(by)}{L(y)}\underset{y\rightarrow +\infty}{\longrightarrow}1, \quad \frac{\tilde{L}(by)}{\tilde{L}(y)}\underset{y\rightarrow +\infty}{\longrightarrow}1.
$$
\end{rem}

\subsection{Generating equations on a plane}\label{section Generalized CIR equations on a plane}

In this section we characterize all equations \eqref{rownanie 2}, with $d=2$, which generate affine models by 
a direct description of the classes $\mathbb{A}_1(a,b; \alpha_1;\eta_1)$ and $\mathbb{A}_2(a,b;\alpha_1,\alpha_2; \eta_1,\eta_2)$.
Our analysis requires an additional regularity assumption that the components of $G$ are strictly positive outside zero and 
\begin{gather}\label{iloraz G regularny}
	\frac{G_2(\cdot)}{G_1(\cdot)}\in C^1(0,+\infty).
\end{gather}
Then $\mathbb{A}_1(a,b; \alpha_1;\eta_1)$ consists of the following equations\\
$$
\bullet \quad \dd R(t)=(aR(t)+b)\dd t+c_0 R(t)^{1/\alpha_1}\Big(G_1 \dd Z_1(t)+G_2 \dd Z_2(t)\Big), 
$$
where $c_0=(\frac{\eta_1}{c_{\alpha_1}})^{\frac{1}{\alpha_1}}$, $G_1,G_2$ are positive constants and $G_1 Z_1(t)+G_2 Z_2(t)$ is an $\alpha_1$-stable process,
$$
\bullet \quad \dd R(t)=(aR(t)+b)\dd t+G_1(R(t-))\dd Z_1(t)+\left(\frac{\eta_1 R(t-)-c_1 G_1^{\alpha_1}(R(t-))}{c_2}\right)^{1/\alpha_1}\dd Z_2(t), 
$$
where $c_1,c_2>0$, $G_1(\cdot)$ is any function such that  
$$
G_1(x)> 0, \quad \frac{\eta_1 x-c_1 G_1^{\alpha_1}(x)}{c_2}>0, \qquad x>0,
$$
and $Z_1, Z_2$ are stable processes with index $\alpha_1$.\\
\noindent
The class $\mathbb{A}_2(a,b;\alpha_1,\alpha_2; \eta_1,\eta_2)$ is a singleton.

The classification above follows directly from the following result.

\begin{tw}\label{tw d=2 independent coord.}
Let $G(x)=(G_1(x), G_2(x))$ be continuous functions such that $G_1(x)>0,G_2(x)>0, x>0$ and \eqref{iloraz G regularny} holds.
Let $Z(t)=(Z_1(t),Z_2(t))$ have independent coordinates of infinite variation with Laplace exponents varying regularly at zero
with indices $\alpha_1,\alpha_2$, respectively, where $2\geq \alpha_1\geq\alpha_2>1$. 
\begin{enumerate}[I)]
\item If $\tilde{J}_{\mu}$  is of the form 
\begin{gather}\label{pierwszy przyp Jmu}
	\tilde{J}_{\mu}(b)=\eta_1 b^{\alpha_1}, \quad b\geq 0,
\end{gather}
with $\eta_1>0, 1<\alpha_1\leq 2$, then $(G,Z)$ is a  generating pair if and only if one of the following two cases holds:
\begin{enumerate}[a)]
\item 
\begin{gather}\label{wspolinionwosc G}
G(x)=c_0 \ x^{1/\alpha_1}\cdot\left(
\begin{array}{ccc}
 G_1\\
 G_2, 
\end{array}
\right), \quad x\geq 0,
\end{gather}
where $c_0=(\frac{\eta_1}{c_{\alpha_1}})^{\frac{1}{\alpha_1}}, G_1>0, G_2>0$ and the process
$$
G_1 Z_1(t)+G_2 Z_2(t), \quad t\geq 0,
$$
is $\alpha_1$-stable.
\item $G(x)$ is such that
\begin{gather}\label{Ib}
c_1 G^{\alpha_1}_1(x)+c_2 G^{\alpha_1}_2(x)=\eta_1 x, \quad x\geq 0,
\end{gather}
with some constants $c_1,c_2>0$, and $Z_1,Z_2$ are $\alpha_1$-stable processes.
\end{enumerate}

\item If $\tilde{J}_{\mu}$  is of the form 
\begin{gather}\label{drugi przyp Jmu}
	\tilde{J}_{\mu}(b)=\eta_1 b^{\alpha_1}+\eta_2 b^{\alpha_2},\quad b\geq 0,
\end{gather}

with $\eta_1,\eta_2>0, 2\geq \alpha_1>\alpha_2>1$ then $(G,Z)$ is a  generating pair if and only if
\begin{gather}\label{postac g w drugim prrzypadku}
G_1(x)=\left(\frac{\eta_1}{c_1}  x\right)^{1 / \alpha_1}, \quad G_2(x)=\left(\frac{\eta_2}{d_2}  x\right)^{1 / \alpha_2}, \quad x\geq 0,
\end{gather}
with some $c_1,d_2>0$ and $Z_1$ is $\alpha_1$-stable, $Z_2$ is $\alpha_2$-stable.
\end{enumerate}
\end{tw}

\noindent
{\bf Proof:} In view of Theorem \ref{TwNiez} the generating pairs $(G,Z)$ are such that 
\begin{gather}\label{równanie z Filipovica d=2}
	J_1(bG_1(x))+J_2(bG_2(x))=x \tilde{J}_{\mu}(b), \quad b,x\geq 0,
\end{gather}
where $\tilde{J}_{\mu}$ takes the form \eqref{pierwszy przyp Jmu} or \eqref{drugi przyp Jmu}. We deduce from \eqref{równanie z Filipovica d=2} the form of $G$ and characterize the noise $Z$. First let us consider the case when 
\begin{gather}\label{znikanie pochodnej}
\left(\frac{G_2(x)}{G_1(x)}\right)^\prime=0, \qquad x>0.
\end{gather}
Then 
$G(x)$ can be written in the form
\begin{gather*}
G(x)=g(x)\cdot\left(
\begin{array}{ccc}
 G_1\\
 G_2, 
\end{array}
\right), \quad x\geq 0,
\end{gather*}
with some function $g(x)\geq 0, x\geq 0$, and constants $G_1>0,G_2>0$. Equation \eqref{rownanie 2} 
amounts then to
\begin{align*}
dR(t)&=F(R(t))+g(R(t-)) \left(G_1dZ_1(t)+G_2 dZ_2(t)\right)\\[1ex]
&=F(R(t))+g(R(t-)) d\tilde{Z}(t), \quad t\geq 0,
\end{align*}
which is an equation driven by the one dimensional L\'evy process $\tilde{Z}(t):=G_1 Z_1(t)+G_2 Z_2(t)$. 
It follows that $\tilde{Z}$ is $\alpha_1$-stable with $\alpha_1\in(1,2]$ and that $g(x)=c_0x^{1/ \alpha_1}, c_0>0$. 
Notice that $Z^{G(x)}(t)=c_0 x^{\frac{1}{\alpha_1}}\tilde{Z}$, so $J_{Z^{G(x)}}(b)=c_{\alpha_1}(c_0 x^{\frac{1}{\alpha_1}}b)^{\alpha_1}=x c_0^{\alpha_1}c_{\alpha_1} b^{\alpha_1}$ and $c_0=(\frac{\eta_1}{c^{\alpha_1}})^{\frac{1}{\alpha_1}}$.
Hence \eqref{pierwszy przyp Jmu} holds and this proves $(Ia)$.

If \eqref{znikanie pochodnej} is not satisfied, then
\begin{gather}\label{niezerowanie pochodnej}
\left(\frac{G_2(x)}{G_1(x)}\right)^\prime\neq 0, \quad x\in (\underline{x},\bar{x}),
\end{gather}
for some interval $(\underline{x},\bar{x})\subset (0,+\infty)$. In the rest of the proof we consider this case and
prove $(Ib)$ and $(II)$.

$(Ib)$ From the equation
\begin{gather}\label{rrr}
J_1(bG_1(x))+J_2(bG_2(x))=x\eta_1 b^{\alpha_1}, \quad b\geq 0, \ x\geq 0,
\end{gather}
we explicitly determine unknown functions. Inserting $b/G_1(x)$ for $b$ yields
\begin{gather}\label{rowwww do eliminacji J1}
J_1(b)+J_2\left(b\frac{G_2(x)}{G_1(x)}\right)=\eta_1\frac{x}{G_1^{\alpha_1}(x)}b^{\alpha_1}, \quad b\geq 0, \quad x>0.
\end{gather}
Differentiation over $x$ yields
 $$
J_2^\prime\left(b\frac{G_2(x)}{G_1(x)}\right)\cdot b \left(\frac{G_2(x)}{G_1(x)}\right)^\prime
=\eta_1\left(\frac{x}{G_1^{\alpha_1}(x)}\right)^\prime b^{\alpha_1}, \quad b\geq 0 ,\quad  x>0.
 $$
Using \eqref{niezerowanie pochodnej} and dividing by $\left(\frac{G_2(x)}{G_1(x)}\right)^\prime$ leads to
 $$
J_2^\prime\left(b\frac{G_2(x)}{G_1(x)}\right)\cdot b=\eta_1 \frac{\left(\frac{x}{G^{\alpha_1}_1(x)}\right)^\prime}{\left(\frac{G_2(x)}{G_1(x)}\right)^\prime}\cdot b^{\alpha_1},\quad b\geq 0 ,\quad x\in(\underline{x},\bar{x}).
$$
By inserting $b\frac{G_1(x)}{G_2(x)}$ for $b$ one computes the derivative of $J_2$: 
$$
J_2^\prime(b)=\eta_1 \frac{\left(\frac{x}{G^{\alpha_1}_1(x)}\right)^\prime\left(\frac{G_1(x)}{G_2(x)}\right)^{\alpha_1-1}}{\left(\frac{G_2(x)}{G_1(x)}\right)^\prime}\cdot b^{\alpha_1-1},\quad b>0 , \quad x\in(\underline{x},\bar{x}).
$$
Fixing $x$ and integrating over $b$ provides
\begin{gather}\label{J_2 wylioczona}
J_2(b)=c_2 b^{\alpha_1}, \quad b>0,
\end{gather}
with some $c_2\geq 0$. Actually $c_2>0$ as $Z_2$ is of infinite variation and $J_2$ can not disappear.

By the symmetry of \eqref{rrr} the same conclusion holds for $J_1$, i.e.
\begin{gather}\label{J_1 wylioczona}
J_1(b)=c_1 b^{\alpha_1}, \quad b>0,
\end{gather}
with $c_1>0$. Using \eqref{J_2 wylioczona} and \eqref{J_1 wylioczona} in \eqref{rrr} gives us \eqref{Ib}. This proves $(Ib)$.

 $II)$ Solving the equation 
\begin{gather}\label{rrrrr}
J_1(bG_1(x))+J_2(bG_2(x))=x(\eta_1 b^{\alpha_1}+\eta_2 b^{\alpha_2}), \quad b,x\geq 0,
\end{gather}
in the same way as we solved \eqref{rrr} yields that 
 \begin{gather}\label{J1,J2 podwojne stabilne}
 J_1(b)=c_1 b^{\alpha_1}+c_2 b^{\alpha_2}, \quad  J_2(b)=d_1 b^{\alpha_1}+d_2 b^{\alpha_2}, \quad b\geq 0,
 \end{gather}
 with $c_1,c_2,d_1,d_2\geq 0$, $c_1+c_2>0, d_1+d_2>0$. From \eqref{rrrrr} and \eqref{J1,J2 podwojne stabilne} we can specify the following conditions for $G$:
  \begin{align}\label{aaaaaa}
c_1G_1^{\alpha_1}(x)+d_1G_2^{\alpha_1}(x)&=\eta_1 x,\\[1ex]\label{bbbbbb}
c_2G_1^{\alpha_2}(x)+d_2G_2^{\alpha_2}(x)&=\eta_2 x.
\end{align}
 We will show that $c_1>0, c_2=0, d_1=0, d_2>0$ by excluding the opposite cases.

 If $c_1>0,c_2>0$, one computes from \eqref{aaaaaa}-\eqref{bbbbbb} that
\begin{gather}\label{G2 na dwa sposoby}
G_1(x)=\left(\frac{1}{c_1}(\eta_1x-d_1G_2^{\alpha_1}(x))\right)^{\frac{1}{\alpha_1}}=\left(\frac{1}{c_2}(\eta_2x-d_2G_2^{\alpha_2}(x))\right)^{\frac{1}{\alpha_2}}, \quad x\geq 0.
\end{gather}
 This means that, for each $x\geq 0$, the value $G_2(x)$ is a solution of the following equation of the $y$-variable
\begin{gather}\label{rownanie z y}
\left(\frac{1}{c_1}(\eta_1x-d_1y^{\alpha_1})\right)^{\frac{1}{\alpha_1}}=\left(\frac{1}{c_2}(\eta_2x-d_2y^{\alpha_2})\right)^{\frac{1}{\alpha_2}}, 
\end{gather} 
with $y\in \left[0,\left(\frac{\gamma_1 x}{d_1}\right)^{\frac{1}{\alpha_1}}\wedge
\left(\frac{\gamma_2 x}{d_2}\right)^{\frac{1}{\alpha_2}}\right]$.
If $d_1=0$ or $d_2=0$ we compute $y=y(x)$ from \eqref{rownanie z y} and see that $d_1y^{\alpha_1}$ or $d_2y^{\alpha_2}$ must be negative either for $x$ sufficiently close to $0$ or $x$ sufficiently large. Now we need to exclude the case $d_1>0,d_2>0$. However,  in the case $c_1, c_2,d_1,d_2>0$ equation \eqref{rownanie z y} has no solutions because, for sufficiently large $x>0$, the left side of \eqref{rownanie z y} is strictly less then the right side. This inequality follows from Proposition \ref{prop o braku rozwiazan} proven below. 

So, we proved that $c_1\cdot c_2=0$ and similarly one proves that $d_1\cdot d_2=0$.
The case $c_1=0, c_2>0, d_1>0, d_2=0$ can be rejected because then $J_1$ would vary regularly with index $\alpha_2$
and $J_2$ with index $\alpha_1$, which is a contradiction. It follows that $c_1>0, c_2=0, d_1=0, d_2>0$ and in this case we obtain \eqref{postac g w drugim prrzypadku} from \eqref{aaaaaa} and \eqref{bbbbbb}. \hfill$\square$
\begin{prop}\label{prop o braku rozwiazan}
Let  $a,b,c,d>0$, $\gamma\in(0,1)$,  $2\geq \alpha_1>\alpha_2>1$. Then for sufficiently large $x>0$ the following inequalities are true
\begin{gather}\label{niernier pomocnicza}
\Big(ax-(bx-cz)^\gamma\Big)^{\frac{1}{\gamma}}-dz> 0, \qquad z\in \Big[0,\frac{b}{c}x\Big],
\end{gather}
\begin{gather}\label{niernier pomocnicza wlasciwa}
(bx-cy^{\alpha_1})^{\frac{1}{\alpha_1}}<(ax-dy^{\alpha_2})^{\frac{1}{\alpha_2}}, \quad y\in \Big[0,\Big(\frac{b}{c}x\Big)^{\frac{1}{\alpha_1}}\wedge \Big(\frac{a}{d}x\Big)^{\frac{1}{\alpha_2}}\Big].
\end{gather}
\end{prop}
{\bf Proof:} First we prove \eqref{niernier pomocnicza} and write it in the equivalent form
\begin{gather}\label{equivalentt}
ax\geq (dz)^{\gamma}+(bx-cz)^\gamma=:h(z).
\end{gather}
Since
$$
h^{\prime}(z)=\gamma \Big(d^\gamma z^{\gamma-1}-c(bx-cz)^{\gamma-1}\Big),
$$
$$
h^{\prime\prime}(z)=\gamma(\gamma-1)\Big(d^\gamma z^{\gamma-2}+c^2(bx-cz)^{\gamma-2}\Big)< 0, \quad z\in \Big[0,\frac{b}{c}x\Big],
$$
the function $h$ is concave and attains its maximum at point
$$
z_0:=\theta x:=\frac{b c^{\frac{1}{\gamma-1}}}{d^{\frac{\gamma}{\gamma-1}}+c^{\frac{\gamma}{\gamma-1}}}x \in \Big[0,\frac{b}{c}x\Big],
$$
which is a root of $h^\prime$. It follows that
\begin{align*}
h(z)\leq h(\theta x)&=(\theta x)^{\gamma}+(bx-c\theta x)^{\gamma}\\
&=(\theta^\gamma+(b-c\theta)^\gamma)x^{\gamma}< ax,
\end{align*}
provided that $x$ is sufficiently large and \eqref{niernier pomocnicza} follows. \eqref{niernier pomocnicza wlasciwa} follows from
\eqref{niernier pomocnicza} by setting $\gamma=\alpha_2/\alpha_1$, $z=y^{\alpha_1}$.
\hfill$\square$

\subsection{An example in 3D}\label{section Example in higher dimensions}

In Section \ref{section Generalized CIR equations on a plane} we proved that in the case $d=2$ the set $\mathbb{A}_2(a,b;\alpha_1,\alpha_2;\eta_1,\eta_2)$ is a singleton. Here we show that  
this property breaks down when $d=3$. In the example below we construct a family of generating pairs 
$(G,Z)$ such that 
\begin{gather}
	J_{Z^{G(x)}}(b)=x\left(\eta_1 b^{\alpha_1}+\eta_2 b^{\alpha_2}\right),\quad b\geq 0,
\end{gather}
with $\eta_1,\eta_2>0, 2\geq \alpha_1>\alpha_2>1$ and such that the related generating equations differ from the canonical representation of $\mathbb{A}_2(a,b;\alpha_1,\alpha_2;\eta_1,\eta_2)$.

\begin{ex} 
Let us consider a process $Z(t)=(Z_1(t),Z_2(t),Z_3(t))$ with independent coordinates such that  $Z_1$ is $\alpha_1$-stable, $Z_2$ is $\alpha_2$-stable,  $Z_3$ is a sum of an $\alpha_1$-  and $\alpha_2$-stable processes. Then
$$
J_1(b)=\gamma_1 b^{\alpha_1}, \quad J_2(b)=\gamma_2 b^{\alpha_2}, \quad J_3(b)=\gamma_3b^{\alpha_1}+\tilde{\gamma}_3 b^{\alpha_2}, \quad b\geq 0,
$$
where $\gamma_1>0,\gamma_2>0,\gamma_3>0,\tilde{\gamma}_3>0$. We are looking for non-negative functions $G_1, G_2,G_3$ solving the equation
\begin{gather}\label{krokodyl}
J_1(bG_1(x))+J_2(bG_2(x))+J_3(bG_3(x))=x \left(\eta_1 b^{\alpha_1}+\eta_2 b^{\alpha_2}\right), \quad x,b\geq 0.
\end{gather}
It follows from \eqref{krokodyl} that
$$
\gamma_1b^{\alpha_1}(G_1(x))^{\alpha_1}+\gamma_2b^{\alpha_2}(G_2(x))^{\alpha_2}
+\gamma_3b^{\alpha_1}(G_3(x))^{\alpha_1}+\tilde{\gamma}_3 b^{\alpha_2}(G_3(x))^{\alpha_2}=
x\left[\eta_1 b^{\alpha_1}+\eta_2b^{\alpha_2}\right], \quad x,b\geq 0,
$$
and, consequently,
$$
b^{\alpha_1}\left[\gamma_1G_1^{\alpha_1}(x)+\gamma_3 G_3^{\alpha_1}(x)\right]+
b^{\alpha_2}\left[\gamma_2G_2^{\alpha_2}(x)+\tilde{\gamma}_3 G_3^{\alpha_2}(x)\right]=x\left[\eta_1 b^{\alpha_1}+\eta_2b^{\alpha_2}\right], \quad x,b\geq 0.
$$
Thus we obtain the following system of equations 
\begin{gather*}
\gamma_1G_1^{\alpha_1}(x)+\gamma_3 G_3^{\alpha_1}(x)=x\eta_1, \\
\gamma_2G_2^{\alpha_2}(x)+\tilde{\gamma}_3 G_3^{\alpha_2}(x)=x\eta_2,
\end{gather*}
which allows us to determine $G_1$ and $G_2$ in terms of $G_3$, that is
\begin{gather}\label{G1}
G_1(x)=\left(\frac{1}{\gamma_1}\left(x\eta_1-\gamma_3 G_3^{\alpha_1}(x)\right)\right)^{\frac{1}{\alpha_1}}\\\label{G2}
G_2(x)=\left(\frac{1}{\gamma_2}\left(x\eta_2-\tilde{\gamma}_3 G_3^{\alpha_2}(x)\right)\right)^{\frac{1}{\alpha_2}}.
\end{gather}
The positivity of $G_1,G_2,G_3$  means that $G_3$ satisfies
\begin{gather}\label{war nieuj}
0\leq G_3(x)\leq \left(\frac{\eta_1}{\gamma_3}x\right)^{\frac{1}{\alpha_1}}\wedge \left(\frac{\eta_2}{\tilde{\gamma}_3}x\right)^{\frac{1}{\alpha_2}}, \quad x\geq 0.
\end{gather}
It follows that $(G,Z)$ with any $G_3$ satisfying \eqref{war nieuj} and $G_1,G_2$ given by \eqref{G1}, \eqref{G2} constitutes a generating pair.
\end{ex}

\section{Applications}\label{section Applications}

Motivated by the form of canonical representations \eqref{rownanie sklajane}  we focus now on the equation
\begin{gather}\label{rownanie numeryka}
dR(t)=(aR(t)+b)\dd t+\sum_{i=1}^{g}d_i^{1/\alpha_i}R(t-)^{1/\alpha_i}\dd Z^{\alpha_i}(t), \quad R(0)=R_0, \ t>0,
\end{gather}
where $a\in\mathbb{R}, b\geq 0, d_i>0$ and $Z^{\alpha_i}$ is an $\alpha_i$-stable process with 
$2\geq \alpha_1>\alpha_2>...>\alpha_g>1$ and $g\geq 1$. By Proposition \ref{prop canonical representation}, \eqref{rownanie numeryka} is the canonical representation of the class $\mathbb{A}_{g}(a,b;\alpha_1,...,\alpha_g;\eta_1,...,\eta_g)$ where 
\begin{gather}\label{ety w kanonicznym}
\eta_i:=c_{\alpha_i}\cdot d_i, 
\end{gather}
and $c_{\alpha_i}$ is given by \eqref{c alpha}.  After characterizing bond prices in the resulted affine model we investigate the flexibility of fitting of \eqref{rownanie numeryka} to risk-free market curves. Our numerical implementations show better performance of \eqref{rownanie numeryka} 
in comparison to the standard CIR equation \eqref{CIR equation}.

Let us start with recalling the concept of pricing based on the semigroup
\begin{gather}\label{semigroup}
\mathcal{Q}_tf(x):= \mathbb{E}[e^{-\int_{0}^{t}R(s)ds}f(R(t))\mid R(0)=x], \quad t\geq 0,
\end{gather}
which was developed in \cite{FilipovicATS}. The formula provides the price at time $0$ of the claim $f(R(t))$ paid at time $t$ given $R(0)=x$. By Theorem 5.3 in \cite{FilipovicATS} for $f_{\lambda}(x):=e^{-\lambda x}, \lambda \geq 0$ we know that
\begin{gather}\label{wycena ogolnie}
\mathcal{Q}_t f_\lambda (x)=e^{-\rho(t,\lambda)-\sigma(t,\lambda)x}, \quad x\geq 0,
\end{gather}
where $\sigma(\cdot,\cdot)$ satisfies the equation
$$
\frac{\partial\sigma}{\partial t}(t,\lambda)=\mathcal{R}(\sigma(t,\lambda)), \quad \sigma(0,\lambda)=\lambda,
$$
and $\rho(\cdot,\cdot)$ is given by
$$
\rho(t,\lambda)=\int_{0}^{t}\mathcal{F}(\sigma(s,\lambda))ds.
$$
The functions $\mathcal{R}, \mathcal{F}$ depend on the generator of $R$, which for \eqref{rownanie numeryka}
takes the form
\begin{align*}
\mathcal{A}f(x)=cx f^{\prime\prime}(x)&+\Big[x\Big(a+\int_{(1,+\infty)}(1 -v)x\tilde{\mu}(\dd v)\Big)+b\Big]f^{\prime}(x)\\[1ex]
&+\int_{(0,+\infty)}[f(x+v)-f(x)-f^{\prime}(x)(1\wedge v)]x\tilde{\mu}(\dd v),
\end{align*}
where
\begin{gather}\label{postac muuu}
\tilde{\mu}(\dd v):= \frac{d_l}{v^{1+\alpha_l}}dv+...+\frac{d_g}{v^{1+\alpha_g}}dv, \quad v>0.
\end{gather}
Recall, if $\alpha_1=2$, then $c=d_1/2$ and $l=2$. Otherwise $c=0$ and $l=1$. Then
\begin{align}\label{wzor na F}\nonumber
\mathcal{R}(\lambda)&:=-c\lambda^2+\Big[a+\int_{(1,+\infty)}(1-v)\tilde{\mu}(\dd v)\Big]\lambda+1+\int_{0}^{+\infty}(1-e^{-\lambda v}-\lambda(1\wedge v))\tilde{\mu}(\dd v),\\[1ex]
\mathcal{F}(\lambda)&:=b\lambda.
\end{align} 
Using \eqref{postac muuu} yields
 \begin{align}\label{wzor na R}\nonumber
\mathcal{R}(\lambda)&=-c\lambda^2+\Big[a+\int_{(1,+\infty)}(1-v)\tilde{\mu}(\dd v)\Big]\lambda+1-\int_{0}^{+\infty}(e^{-\lambda v}-1+\lambda v)\tilde{\mu}(\dd v)\\[1ex]\nonumber
&\quad -\lambda\int_{(1,+\infty)}(1-v)\tilde{\mu}(\dd v)=-c\lambda^2+a\lambda+1-\sum_{i=l}^g\eta_k\lambda^{\alpha_k}\\[1ex]
&=1+a\lambda -\sum_{i=1}^g\eta_k\lambda^{\alpha_k}.
\end{align}

Application of the pricing procedure above for $f_\lambda$ with $\lambda=0$  
allows us to obtain from \eqref{wycena ogolnie} the prices of zero-coupon bonds. Using the closed form formula 
\eqref{wzor na R} leads to the following result.
\begin{tw}\label{tw o cenach obligacji}
The zero-coupon bond prices in the affine model generated by \eqref{rownanie numeryka} are equal
\begin{gather}\label{affine prices tw}
	P(t,T)=e^{-A(T-t)-B(T-t)R(t)},
\end{gather}
where $B$ and $A$ are such that 
	\begin{align}\label{B indep. coord.}
	B^\prime(v)&=1+aB(v)-\sum_{i=1}^{g}\eta_iB^{\alpha_i}(v),\quad B(0)=0,\\[1ex]\label{A indep. coord.}
	A^\prime(v)&=b B(v), \quad A(0)=0,
	\end{align}
with $\{\eta_i\}$ given by \eqref{ety w kanonicznym}.
\end{tw}

In the case when $g=1$ and $\alpha_1=2$ equation \eqref{B indep. coord.} becomes a Riccati equation and its explicit solution
provides bond prices for the classical CIR equation. In the opposite case \eqref{B indep. coord.} can be solved 
by numerical methods which exploit the tractable form of the function $\mathcal{R}$ given by \eqref{wzor na R}. Note that $\mathcal{R}$ is continuous, $\mathcal{R}(0)=1$ and  $\lim_{\lambda\rightarrow +\infty}\mathcal{R}(\lambda)=-\infty$. Thus 
$\lambda_0:=\inf\{\lambda>0:\mathcal{R}(\lambda)=0\}$
is a positive number and 
\begin{gather}\label{zachowanie R w lambda 0}
\mathcal{R}(\lambda_0)=0, \quad \mathcal{R}^{\prime}(\lambda_0)<0. 
\end{gather}
The function
\begin{gather}\label{funkcja G}
\mathcal{G}(x):=\int_{0}^{x}\frac{1}{\mathcal{R}(y)}dy, \quad x\in[0,\lambda_0),
\end{gather}
is strictly increasing and its behaviour near $\lambda_0$ can be estimated by
substituting $z=\frac{1}{\lambda_0-y}$ in \eqref{funkcja G} and using the inequality
$$
(\lambda_0-h)^{\alpha}\geq \lambda_0^{\alpha}-\alpha\lambda_0^{\alpha-1}h, \quad h\in(0, \lambda_0),\quad \alpha\in(1,2) .
$$
For the case when $\alpha_1=2$ this yields for $x \in [0, \lambda_0)$
\begin{align}\label{oszacowanie dla G}\nonumber
\mathcal{G}(x)&=\int_{1/\lambda_0}^{1/(\lambda_0-x)}\frac{1}{\mathcal{R}(\lambda_0-\frac{1}{z})}\cdot\frac{1}{z^2} \dd z\\[1ex]\nonumber
&=\int_{1/\lambda_0}^{1/(\lambda_0-x)}\frac{1}{z^2+a\lambda_0 z^2-az-\eta_1(\lambda_0 z-1)^2-\sum_{i=2}^{g}\eta_i z^2(\lambda_0-\frac{1}{z})^{\alpha_i}} \ \dd z\\[1ex]\nonumber
&\geq \int_{1/\lambda_0}^{1/(\lambda_0-x)}\frac{1}{z^2+a\lambda_0 z^2-az-\eta_1(\lambda_0 z-1)^2-\sum_{i=2}^{g}\eta_i z^2(\lambda_0^{\alpha_i}-\alpha_i\lambda_0^{\alpha_i-1}\frac{1}{z})} \ \dd z\\[1ex]\nonumber
&=\int_{1/\lambda_0}^{1/(\lambda_0-x)}\frac{1}{z^2(1+a\lambda_0-\eta_1\lambda_0^2-\sum_{i=2}^{g}\eta_i\lambda_0^{\alpha_i})+z(2\eta_1\lambda_0-a+\sum_{i=2}^{g}\alpha_i\eta_i\lambda_0^{\alpha_i-1})-\eta_1} \ \dd z\\[1ex]
&=\int_{1/\lambda_0}^{1/(\lambda_0-x)}\frac{1}{\mathcal{R}(\lambda_0)z^2-\mathcal{R}^{\prime}(\lambda_0)z-\eta_1} \ \dd z.\end{align}
It follows from \eqref{oszacowanie dla G} and \eqref{zachowanie R w lambda 0} that
$$
\lim_{x\rightarrow \lambda_0^{-}}\mathcal{G}(x)=+\infty,
$$ 
so $\mathcal{G}$ is invertible and $\mathcal{G}^{-1}$ exists on $[0,+\infty)$. Writing \eqref{B indep. coord.} as
$$
B^\prime(v)=\mathcal{R}(B(v)), \quad B(0)=0,
$$
we see that 
$$
\frac{d}{dv}\mathcal{G}(B(v))=\frac{1}{\mathcal{R}(B(v))}B^\prime(v)=1,
$$
and consequently
$$
\mathcal{G}(B(v))=v, \quad v\geq 0.
$$
Representing $B(\cdot)$ as the inverse of $\mathcal{G}(\cdot)$ enables its numerical computation. Hence, with  $\mathcal{G}^{-1}(\cdot)$ at hand we can derive bond prices, spot rates and swap rates in the model generated by \eqref{rownanie numeryka}. The dependence of $\mathcal{G}^{-1}(\cdot)$ on the parameters
$a,\alpha_1,...,\alpha_g, \eta_1,...,\eta_g$ plays a central role in the problem of fitting the model to real data. In what follows 
we present the results of calibration of \eqref{rownanie numeryka}  to market quotes of spot rates, Libor and swap rates. 
 
\subsection{Calibration of canonical models to market data}\label{section Calibration}

\begin{figure}[htb]
\begin{center}
\includegraphics[height=2.5in,width=3in,angle=0]{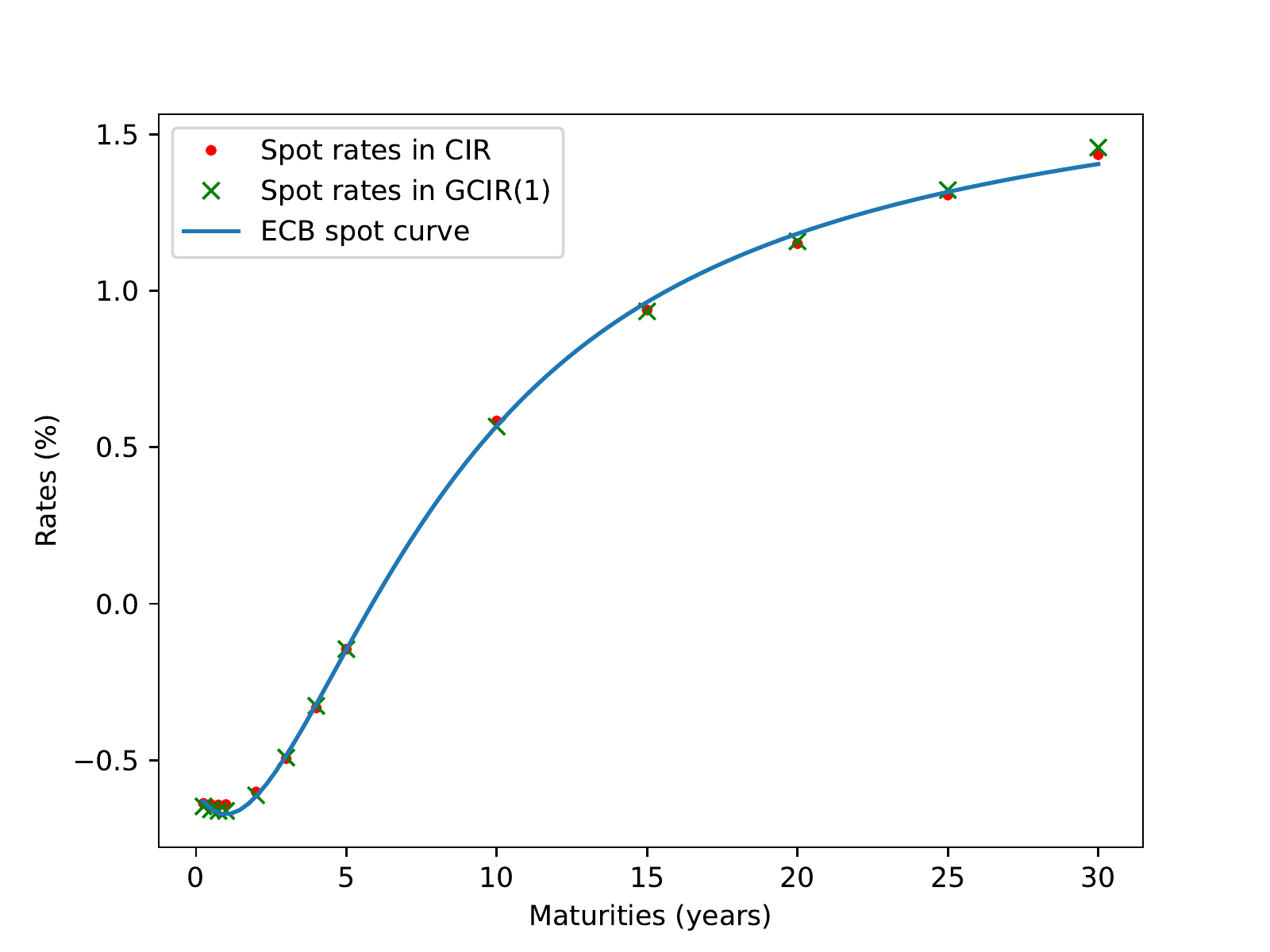}
\includegraphics[height=2.5in,width=3in,angle=0]{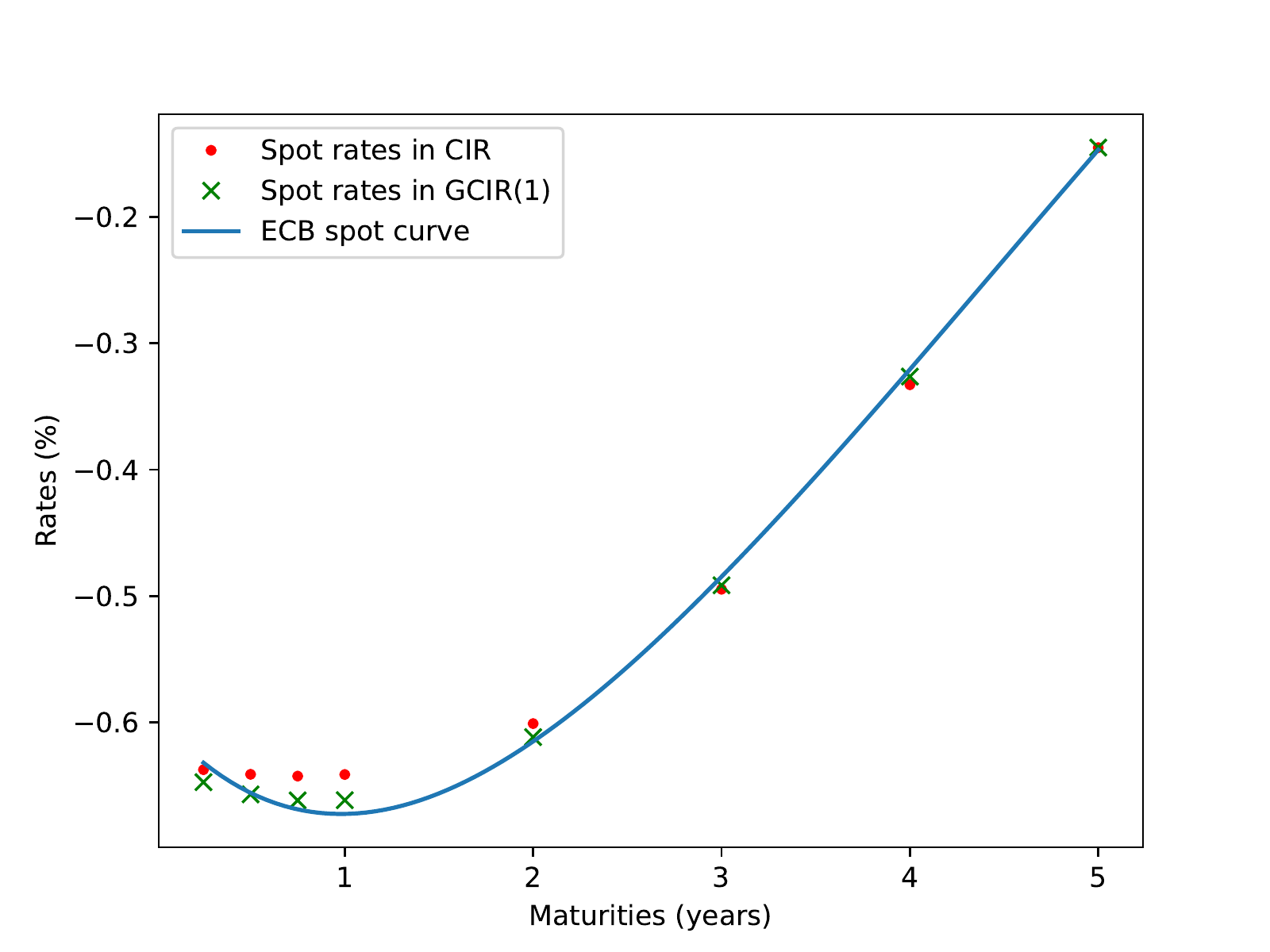}
\caption{Calibration to the ECB curves from 10.01.2018. View for all/small maturities.}\label{plot 2018}
\end{center}
\end{figure}

Our first calibration procedure is concerned with the spot yield curves of European Central Bank (ECB) computed from the zero coupon AAA-rated bonds. The maturity grip consists of $33$ points starting from $3$ months and ending with 30 years. This set was, however, restricted to $13$ points to speed up computations. All maturities less than $5$ years were included to save rapid changes of the curves near zero. A glance at the historical data from 2016 to 2023 reveals significant changes in the shape of curves appearing after March 2022. The classical CIR model could be fitted relatively well to previous curves but performed much worse for the newer ones. In both cases, however, the addition of new stable noise components resulted in reduction of the calibration error.  For a calibration based on maturities $T_1<...<T_M$ the fitting error measures a relative distance of the  model spot rates 
\begin{gather}\label{spot rates in the model}
y(T_i):=\frac{1}{T_i}\left(\frac{1}{P(0,T_i)}-1\right), \quad i=1,2,...,M,
\end{gather}
from the empirical ones $\hat{y}(T_i), i=1,2,...,M$.
It is given by the formula
\begin{gather}\label{error spot rates}
Error(a,b,\alpha_1,...,\alpha_g, d_1,...,d_g):=\sum_{i=1}^{M}\frac{(y(T_i)-\hat{y}(T_i))^2}{\hat{y}^2(T_i)}.
\end{gather}
For the curve from 10.01.2018 we can see that a good fitting of the CIR model can be substantially improved by
replacing the Wiener process by a  stable noise with index $\alpha=1.58$. The effect is strongly apparent especially for small maturities, see Fig. \ref{plot 2018}. The increase of the number of noise components 
causes further decrease of the fitting error but in a lesser extent, see Tab. \ref{table 2018}, where GCIR(g) stands for the generalized CIR equation \eqref{rownanie numeryka} with $g$ components.

\begin{figure}[htb]
\begin{center}
\begin{tabular}{| c | c | l |}   
   \hline
   Model & Calibration error $\times 100$ & Stability indices \\ \hline
    CIR & 0.95141785& $\alpha=2$ \\ \hline
       GCIR(1)& 0.44735953&$\alpha=1.58$\\ \hline
    GCIR(2)& 0.44505444&$\alpha_1=2$, $\alpha_2=1.53$\\ \hline
        GCIR(3)&0.44148324&$\alpha_1=2 $, $\alpha_2=1.91 $, $\alpha_3=1.42$\\ \hline
            GCIR(4)&0.43932515&$\alpha_1=2 $,  $\alpha_2=1.45 $, $\alpha_3=1.44$, $\alpha_4=1.29 $ \\ \hline
                GCIR(5)& 0.43918035&$\alpha_1=2 $, $\alpha_2=1.315 $, $\alpha_3=1.311$, $\alpha_4=1.308$, $\alpha_5=1.23$ \\ \hline
    \end{tabular}
    \captionof{table}{Error reduction - calibration to the ECB rates from 10.01.2018.}\label{table 2018}
\end{center}
\end{figure}

\begin{figure}[htb]
\begin{center}
\includegraphics[height=2.5in,width=3in,angle=0]{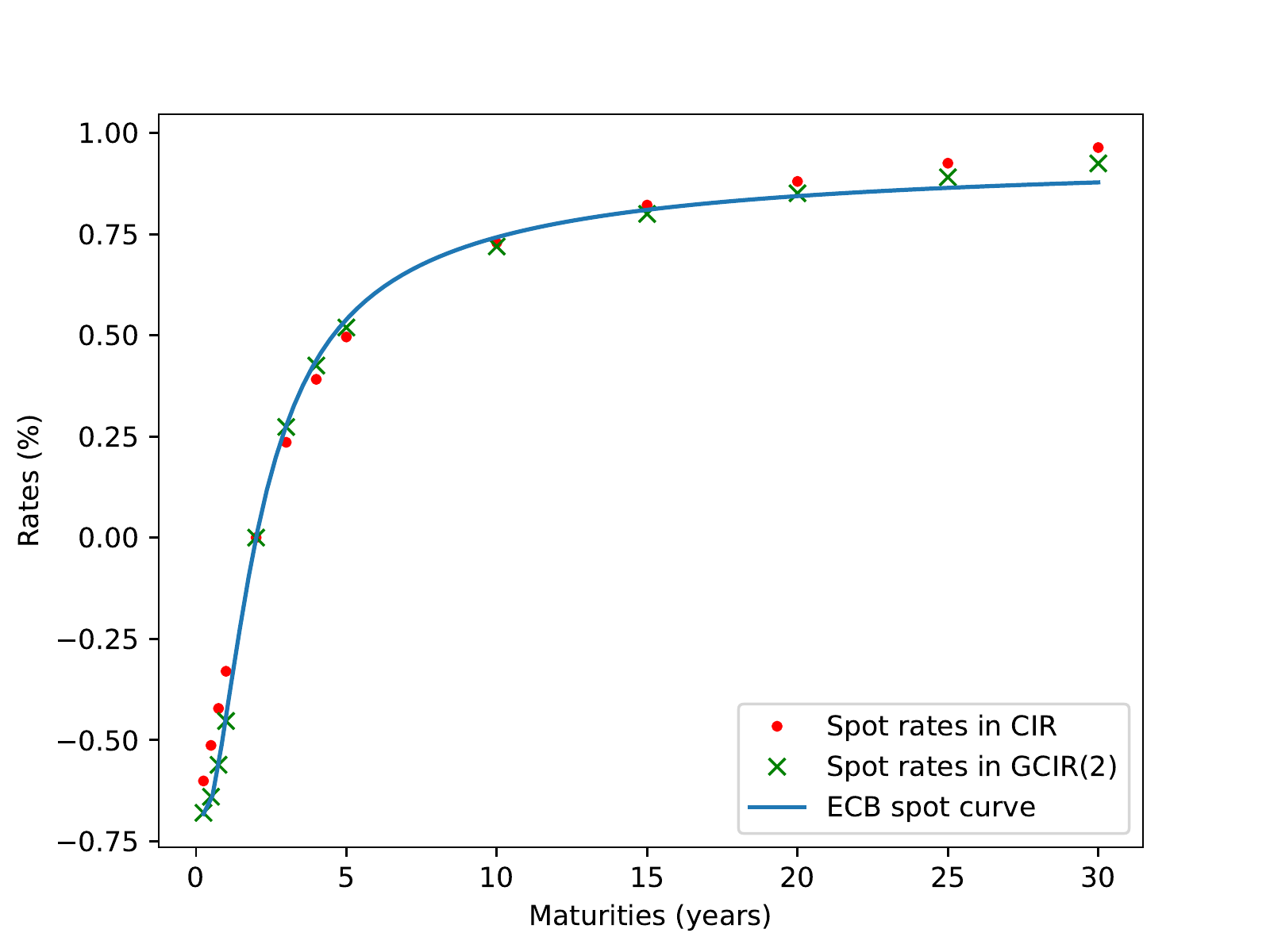}\\
\includegraphics[height=2.5in,width=3in,angle=0]{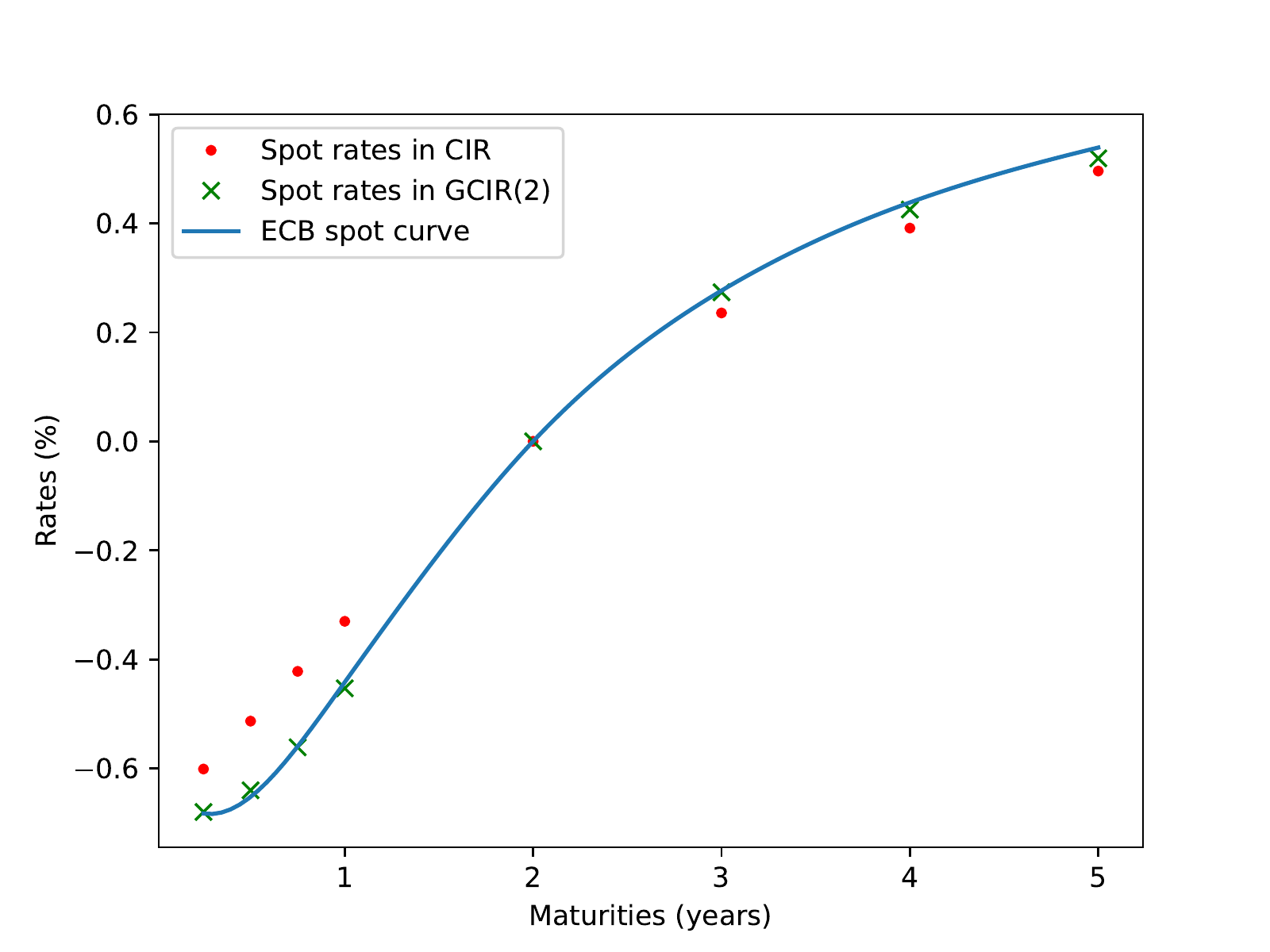}
\includegraphics[height=2.5in,width=3in,angle=0]{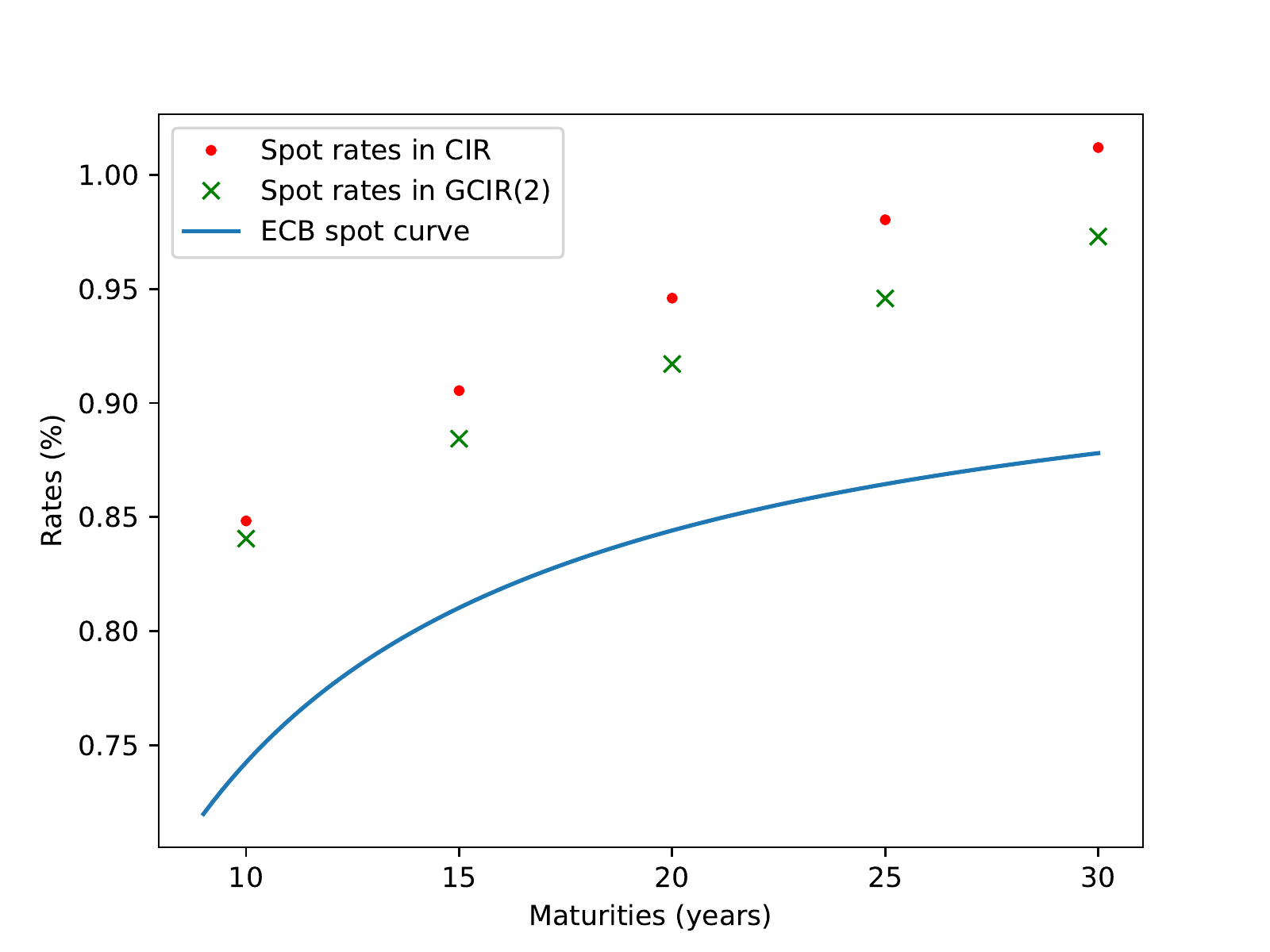}
\caption{Calibration to the ECB curves from 8.04.2022. View for all/small/large maturities.}\label{plot 2022}
\end{center}
\end{figure}

For the data from 8.04.2022 the CIR model turned out to be the most efficient among one dimensional models, though the fitting error is much greater then in the previous example, see Tab. \ref{table 2022} and Fig. \ref{plot 2022}. Models with higher noise dimension provide, however, better results starting from the gratest error reduction by the alpha-CIR model of \cite{JiaoMaScotti} with $\alpha=1.04$.

\begin{figure}[htb]
\begin{center}
    \begin{tabular}{| c | c | l |}
    \hline    
     Model & Calibration error  $\times 100$ & Stability indices \\ \hline
    CIR & 24.10280133& $\alpha=2$ \\ \hline
    GCIR(2)& 0.83059934&$\alpha_1=2$, $\alpha_2=1.99$\\ \hline
    GCIR(3)& 0.83055904&$\alpha_1=2$, $\alpha_2=1.17$, $\alpha_3=1.14$\\ \hline
    GCIR(4)&0.83050323&$\alpha_1=2$, $\alpha_2=1.35$, $\alpha_3=1.25$, $\alpha_4= 1.21$\\ \hline
    GCIR(5)&0.83049801&$\alpha_1=2$, $\alpha_2=1.53$, $\alpha_3=1.48$, $\alpha_4=1.35$, $\alpha_5=1.23$\\ \hline
    \end{tabular}
    \captionof{table}{Error reduction - calibration to the ECB rates from 8.04.2022.}\label{table 2022}    
\end{center}
\end{figure}

\begin{figure}[htb]
\begin{center}
\includegraphics[height=2.5in,width=3in,angle=0]{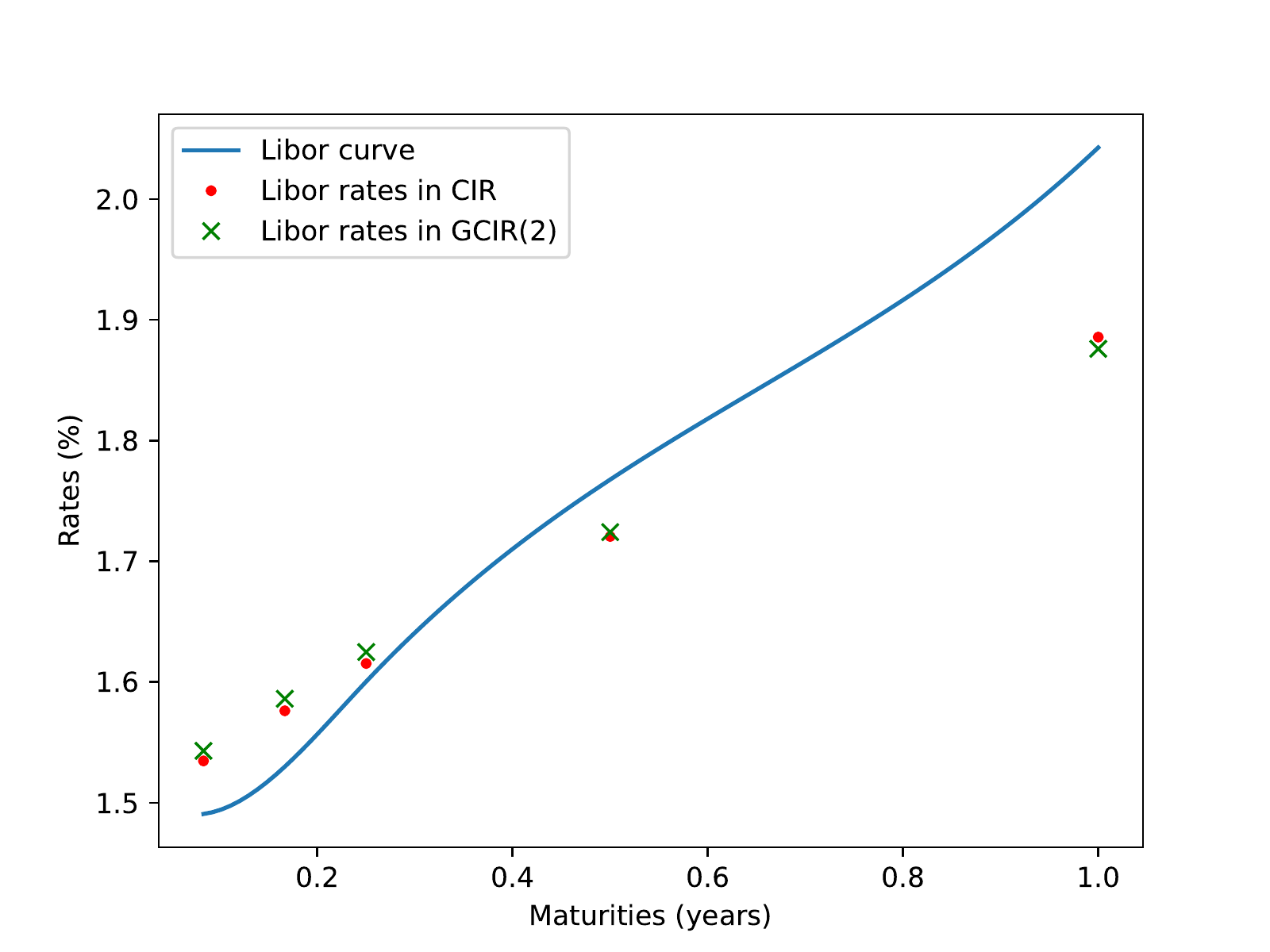}
\includegraphics[height=2.5in,width=3in,angle=0]{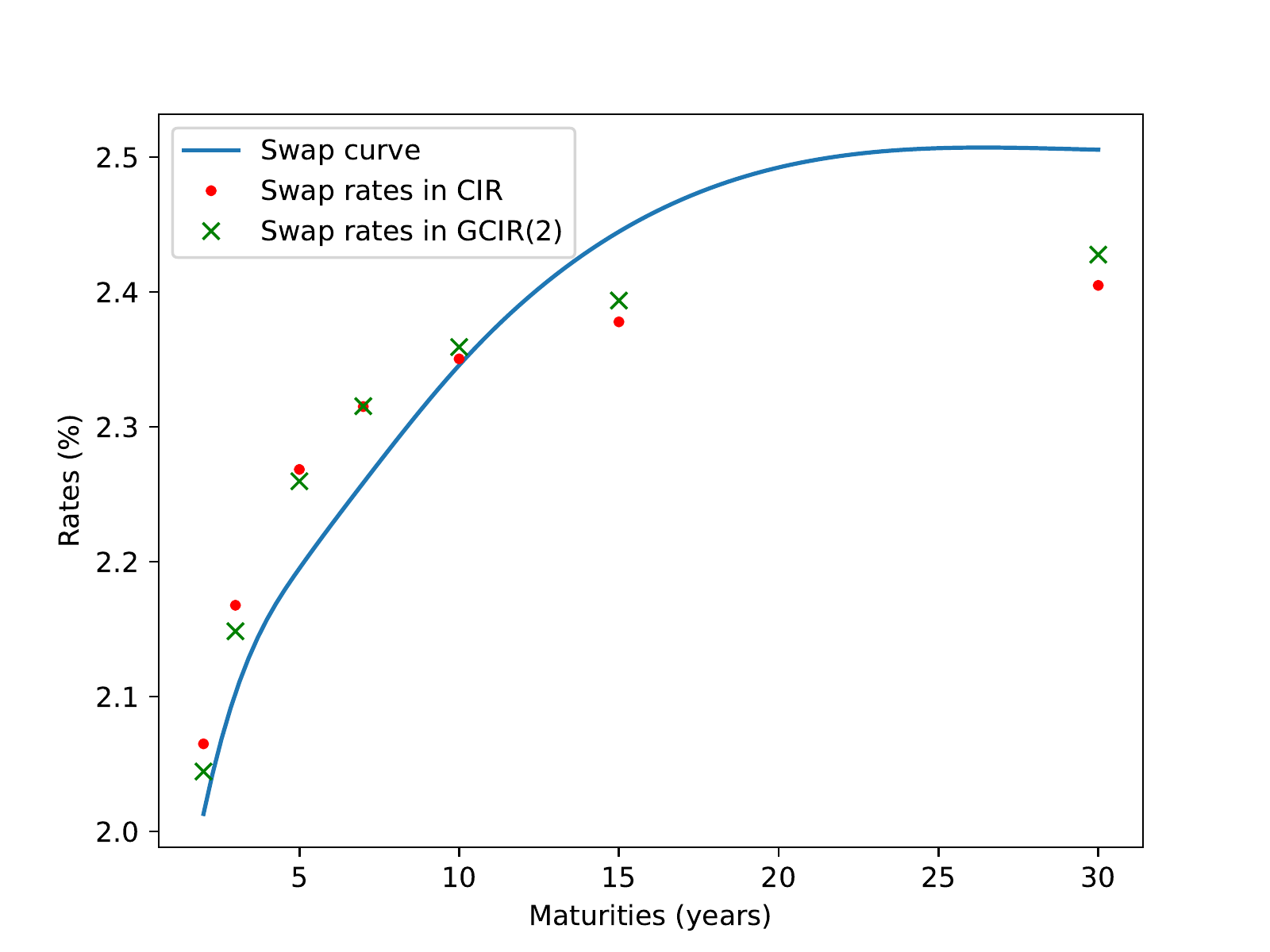}
\caption{Calibration to the Libor and swap curves from 14.12.2017}\label{plot swap libor 2017}
\end{center}
\end{figure}

Our second calibration procedure was based on Libor and $6$-months swap rates with maturities resp. $\{T_i\}, i=1,...,M_1$ and $\{U_i\}, i=1,...,M_2$. The term structure of interest rates for maturities below one year are represented by Libor quotes while swap rates correspond to selected maturities from 1 year up to 30 years. A direct extention of \eqref{error spot rates} 
leads to the calibration error of the form
$$
Error(a,b,\alpha_1,...,\alpha_g, d_1,...,d_g):=\sum_{i=1}^{M_1}\frac{(L(T_i)-\widehat{L}(T_i))^2}{\widehat{L}^2(T_i)}+
\sum_{i=1}^{M_2}\frac{(S(U_i)-\widehat{S}(U_i))^2}{\widehat{S}^2(U_i)},
$$
where Libor rates $L(T_i)$ are defined like \eqref{spot rates in the model} and swap rates by 
$$
S(U_i)=\frac{1-P(0,U_i)}{\frac{1}{2}\sum_{k=1}^{i}P(0,U_k)}, \quad i=1,...,M_2.
$$
The best one dimensional model for the data from 14.12.2017 was CIR, but, again, multivariate models generated better results. The passage from $g=1$ to $g=2$, i.e. to the  $\alpha$-CIR model with $\alpha=1.16$, gave the highest error reduction, which was particularly effective for the swap rates. All of them were pushed closer the empirical swap curve.  The results are presented in Fig.\ref{plot swap libor 2017} and Tab. \ref{table swap libor 2017}.

\begin{figure}[htb]
\begin{center}
    \begin{tabular}{| c  | c | c | c | l |}
    \hline
    Model & \begin{tabular}{c} Calibration error \\$\times 100$
    \end{tabular} 
       &\begin{tabular}{c} Libor error \\$\times 100$  \end{tabular} & \begin{tabular}{c} Swap error \\ $\times 100$  \end{tabular} & Stability indices \\ \hline
    CIR & 1.42225593 &0.84831146&0.57394447& $\alpha=2$ \\ \hline
  GCIR(2)& 1.37316050&1.00280671&0.37035379&$\alpha_1=2$, $\alpha_2=1.16$\\ \hline
  GCIR(3)& 1.37309034&1.00987818&0.36321216&$\alpha_1=2$, $\alpha_2=1.94$, $\alpha_3= 1.15$\\ \hline
  GCIR(4)&1.37308709&1.00989365&0.36319344& \begin{tabular}{c} $\alpha_1=2$, $\alpha_2=1.99$, \\ $\alpha_3= 1.54$, $\alpha_4=1.15$ \end{tabular} \\ \hline
 \end{tabular}
    \captionof{table}{Error reduction - calibration to the Libor and swap rates from 14.12.2017.}\label{table swap libor 2017}    
     \end{center}
\end{figure}

\subsubsection{Remarks on computational methodology}

Our computation were performed in the Python programming language. The calibration error was minimized with the use of the Nelder-Mead algorithm which turned out to be most effective among all available algorithms for local minimization in the Python library. The computation time of calibration which depends, of course, on the number of noise components, lied in the range 100-13.000 seconds but often 
did not exceed 800 seconds. This stays in a strong contrast  to the CIR model for which the closed form formulas 
shorten the calibration to the 2 second limit. We suspect that global optimization algorithms would provide even better fit, but they were to slow for the data with more than several maturities.


\newpage

\section{Appendix}

\noindent{\bf Proof of Proposition \ref{prop wstepny}:}
$(A)$ It was shown in \cite [Theorem 5.3]{FilipovicATS} that the generator of a general positive Markovian short rate generating  an affine model is of the form
\begin{align}\label{generator Filipovica}
\mathcal{A}f(x)=&c x f^{\prime\prime}(x)+(\beta x+\gamma)f^\prime(x)\\[1ex] \nonumber
&+\int_{(0,+\infty)}\Big(f(x+y)-f(x)-f^\prime(x)(1\wedge y)\Big)(m(\dd y)+x\mu(\dd y)), \quad x\geq 0,
\end{align}
for $f\in\mathcal{L}(\Lambda)\cup C_c^2(\mathbb{R}_{+})$, where 
$\mathcal{L}(\Lambda)$ is the linear hull of $\Lambda:=\{f_\lambda:=e^{-\lambda x}, \lambda\in(0,+\infty)\}$
and $C_c^2(\mathbb{R}_{+})$ stands for the set of twice continuously differentiable functions with compact support in $[0,+\infty)$. 
Above $c, \gamma\geq 0$, $\beta\in\mathbb{R}$ and $m(\dd y)$, $\mu(\dd y)$ are nonnegative Borel measures on $(0,+\infty)$ satisfying
\begin{gather}\label{warunki na iary Filipovica}
\int_{(0,+\infty)}(1\wedge y)m(\dd y)+\int_{(0,+\infty)}(1\wedge y^2)\mu(\dd y)<+\infty.
\end{gather}

The generator of the short rate process given by \eqref{rownanie 2} equals
\begin{align*}
 \mathcal{A}_{R}f(x) =& f^\prime(x)F(x)+\frac{1}{2}f^{\prime\prime}(x)\langle QG(x),G(x)\rangle \\
& +\int_{\mathbb{R}^d}\Big(f(x+\langle G(x),y\rangle)-f(x)-f^\prime(x)\langle G(x),y\rangle\Big)\nu(\dd y) \\
 = & f^\prime(x)F(x)+\frac{1}{2}f^{\prime\prime}(x)\langle QG(x),G(x)\rangle \\
&+\int_{\mathbb{R}}\Big(f(x+v)-f(x)-f^\prime(x)v\Big)\nu_{G(x)}(\dd v)
\end{align*}
where  $f$ is a bounded, twice continuously differentiable function. 

By Proposition \ref{prop o skokach z dodatniosci rozwiazania} below, the support of the measure $\nu_{G(x)}$ is contained in $[-x,+\ns)$, thus it follows that  
\begin{align}\label{Generatorr R gen}\nonumber
\mathcal{A}_{R}f(x) = &f^\prime(x)F(x)+\frac{1}{2}f^{\prime\prime}(x)\langle QG(x),G(x)\rangle \\ \nonumber
&+\int_{(0, +\ns)}\Big(f(x+v)-f(x)-f^\prime(x)(1\wedge v) \Big)\nu_{G(x)}(\dd v)\\ \nonumber
&+f^\prime(x)\int_{(0, +\ns)}\Big((1\wedge v)-v\Big)\nu_{G(x)}(\dd v)\\[1ex] \nonumber 
&+\int_{(-\ns, 0)}\Big(f(x+v)-f(x)-f^\prime(x)v \Big)\nu_{G(x)}(\dd v)\\ \nonumber
= & \frac{1}{2}f^{\prime\prime}(x)\langle QG(x),G(x)\rangle + f^\prime(x)\sbr{F(x)+\int_{(1,+\ns)}\Big(1- v\Big)\nu_{G(x)}(\dd v)} \\[1ex] \nonumber
&+\int_{(0, +\ns)}\Big(f(x+v)-f(x)-f^\prime(x)(1\wedge v) \Big)\nu_{G(x)}(\dd v)\\ 
&+\int_{[-x, 0)}\Big(f(x+v)-f(x)-f^\prime(x)v \Big)\nu_{G(x)}(\dd v).  
\end{align}

Comparing \eqref{Generatorr R gen} with \eqref{generator Filipovica} applied to a function $f_{\lambda}$ with $\lambda >0$ such that $f_{\lambda}(x) = e^{-\lambda x}$ for $x \ge 0$, we get

\begin{align}
&c x \lambda^2  - (\beta x+\gamma)  \lambda  \nonumber \\[1ex] \nonumber
&+ \int_{(0,+\infty)}\Big(e^{-\lambda y}-1+\lambda (1\wedge y)\Big)(m(\dd y)+x\mu(\dd y)) \\[1ex] \nonumber
& -  \frac{1}{2} \lambda^2 \langle QG(x),G(x)\rangle + \sbr{F(x)+\int_{(1,+\ns)}\Big(1- v\Big)\nu_{G(x)}(\dd v)} \lambda  
\\[1ex] \nonumber
& - \int_{(0,+\infty)}\Big(e^{-\lambda v}-1+\lambda (1\wedge v)\Big)\nu_{G(x)}(\dd v) \\[1ex] \label{uaua}
& =  \int_{[-x, 0)}\Big(e^{-\lambda v}-1+\lambda  v \Big)\nu_{G(x)}(\dd v) , \quad \lambda >0, x \geq 0.
\end{align}
Comparing the left and the right sides of \eqref{uaua} we see that the left side grows no faster than a quadratic polynomial of $\lambda$ while the right side grows faster that $d e^{\lambda y}$ for some $d, y >0$, unless  the support of the measure $\nu_{G(x)}(\dd v)$ is contained in $[0,+\ns)$. It follows that $\nu_{G(x)}(\dd v)$ is concentrated on $[0,+\infty)$, hence $(a)$ follows, and 

\begin{align}
&c x \lambda^2  - (\beta x+\gamma)  \lambda  \nonumber \\[1ex] \nonumber
& -  \frac{1}{2} \lambda^2 \langle QG(x),G(x)\rangle + \sbr{F(x)+\int_{(1,+\ns)}\Big(1- v\Big)\nu_{G(x)}(\dd v)} \lambda  
\\[1ex] \label{uaua1}
& = \int_{(0,+\infty)}\Big(e^{-\lambda y}-1+\lambda (1\wedge y)\Big)\rbr{\nu_{G(x)}(\dd y) - m(\dd y)- x\mu(\dd y)}, \quad \lambda >0, x \geq 0.
\end{align}
Dividing both sides of the last equality by $\lambda^2$ and using the  estimate  
$$\frac{e^{-\lambda y}-1+\lambda (1\wedge y)}{\lambda^2} \le \rbr{\frac{1}{2} y^2} \wedge \rbr{\frac{e^{-\lambda}-1+\lambda}{\lambda^2}}$$ 
we get that that the left side of \eqref{uaua1} converges to $c x -  \frac{1}{2} \langle QG(x),G(x)\rangle$ as $\lambda\rightarrow +\infty$, while the right side converges to  $0$. This yields \eqref{mult. CIR condition}, i.e. 
\begin{align}\label{W1}
c x=& \frac{1}{2}\langle QG(x),G(x)\rangle,\quad x\geq 0.
\end{align}
Next, fixing $x\geq 0$ and comparing \eqref{Generatorr R gen} with \eqref{generator Filipovica} applied to a function from the domains of both generators and such that $f(x) = f'(x) = f ''(x) =0$ we get 
\[
\int_{(0,+\infty)} f(x+y) (m(\dd y)+x\mu(\dd y)) =  \int_{(0, +\ns)} f(x+v) \nu_{G(x)}(\dd v) 
\]
for any such a function,
which yields 
\begin{gather}\label{rozklad nuG na sume}
\nu_{G(x)}(\dd v)\mid_{(0,+\infty)}=m(\dd v)+x\mu(\dd v), \quad x\geq 0.
\end{gather} 
This implies also
\begin{align}
\label{W2}
\beta x+\gamma=&F(x)+\int_{(1,+\ns)}\Big(1 - v\Big)\nu_{G(x)}(\dd v),\quad x\geq 0.
\end{align}

$(b)$ Setting $x=0$ in \eqref{rozklad nuG na sume} yields 
\begin{gather}\label{nu G0}
\nu_{G(0)}(\dd v)\mid_{(0,+\infty)}=m(\dd v).
\end{gather}
To prove \eqref{nu G0 finite variation}, by \eqref{warunki na iary Filipovica} and \eqref{nu G0}, we need to show that
\begin{gather}\label{cocococ}
\int_{(1,+\infty)}v\nu_{G(0)}(\dd v)<+\infty.
\end{gather}
It is true if $G(0)=0$ and for $G(0)\neq 0$ the following estimate holds
\begin{align*}
\int_{(1,+\infty)}v\nu_{G(0)}(\dd v)&=\int_{\mathbb{R}^d}\langle G(0),y\rangle\mathbf{1}_{[1,+\infty)}(\langle G(0),y\rangle)\nu(\dd y)\\[1ex]
&\leq \mid G(0)\mid\int_{\mathbb{R}^d}\mid y\mid\mathbf{1}_{[1/\mid G(0)\mid,+\infty)}(\mid y\mid)\nu(\dd y),
\end{align*}
and \eqref{cocococ} follows from \eqref{warunek na miare Levyego 2}.

$(c)$   \eqref{rozklad nu G(x)} follows from \eqref{rozklad nuG na sume} and \eqref{nu G0}.
To prove \eqref{war calkowe na mu} we use \eqref{rozklad nu G(x)},  \eqref{nu G0 finite variation}
and the following estimate for $x\geq 0$:
\begin{align*}
\int_{0}^{+\infty}(v^2\wedge v)\nu_{G(x)}(\dd v)&=\int_{\mathbb{R}^d}(\mid\langle G(x),y\rangle\mid^2\wedge\langle G(x),y\rangle)\nu(\dd y)\\[1ex]
&\leq \Big(\mid G(x)\mid^2\vee \mid G(x)\mid\Big)\int_{\mathbb{R}^d}(\mid y\mid^2\wedge \mid y\mid)\nu(\dd y)<+\infty,
\end{align*}
In the last line we used \eqref{warunek na miare Levyego 1} and \eqref{warunek na miare Levyego 2}.

$(d)$ It follows from \eqref{W2} and \eqref{rozklad nu G(x)} that 
\begin{align*}
\beta x+\gamma&=F(x)+\int_{(1,+\infty)}(1 - v)\nu_{G(x)}(\dd v)\\[1ex]
&=F(x)+\int_{(1,+\infty)}(1-v)\nu_{G(0)}(\dd v)+x\int_{(1,+\infty)}(1-v)\mu(\dd v), \quad x\geq 0.
\end{align*}
Consequently, \eqref{linear drift}
follows with
$$
a:=\Big(\beta-\int_{(1,+\infty)}(1-v)\mu(\dd v)\Big), \ b:=\Big(\gamma-\int_{(1,+\infty)}(1-v)\nu_{G(0)}(\dd v)\Big),
$$
and $b\geq \int_{(1,+\infty)}(v-1)\nu_{G(0)}(\dd v)$ because $\gamma\geq 0$.

$(B)$ We use \eqref{W2},  \eqref{linear drift} and \eqref{rozklad nuG na sume} to write \eqref{generator Filipovica} in the form
\begin{align*}
\mathcal{A}f(x)=cx f^{\prime\prime}(x)&+\Big[ax +b+\int_{(1,+\infty)}(1 -v)\nu_{G(x)}(\dd v)\Big]f^{\prime}(x)\\[1ex]
&+\int_{(0,+\infty)}[f(x+v)-f(x)-f^{\prime}(x)(1\wedge v)]\nu_{G(x)}(\dd v)\}.
\end{align*}
In view of \eqref{rozklad nuG na sume} and \eqref{nu G0} we see that
\eqref{generator R w tw} is true.

\begin{prop}\label{prop o skokach z dodatniosci rozwiazania}
Let $G:[0,+\infty)\rightarrow \mathbb{R}^d$ be continuous. If the equation \eqref{rownanie 2} has a non-negative strong solution for any initial condition $R(0)=x\geq 0$, then
\begin{gather}\label{ograniczenia skokow}
\forall x\geq 0 \quad \nu{\{y\in\mathbb{R}^d: x+\langle G(x),y\rangle<0\}}=0.
\end{gather}
In particular, the support of the measure $\nu_{G(x)}(\dd v)$ is contained in $[-x,+\infty)$.
\end{prop}
{\bf Proof:} Let us assume to the contrary, that for some $x\geq 0$
$$
\nu{\{y\in\mathbb{R}^d: x+\langle G(x),y\rangle<0\}}>0.
$$
Then there exists $c>0$ such that
$$
\nu{\{y\in\mathbb{R}^d: x+\langle G(x),y\rangle<-c\}}>0.
$$
Let $A\subseteq \{y\in\mathbb{R}^d: x+\langle G(x),y\rangle<-c\}$ be a Borel set separated from zero. By the continuity of $G$ we have that for some $\varepsilon>0$:
\begin{gather}\label{pani w szpileczkach}
\tilde{x}+\langle G(\tilde{x}),y\rangle<-\frac{c}{2}, \quad  \tilde{x}\in[(x-\varepsilon)\vee 0,x+\varepsilon],\quad y\in A.
\end{gather}
Let $Z^2$ be a L\'evy processes with characteristics $(0,0,\nu^2(dy))$, where $\nu^2(dy):=\mathbf{1}_{A}(y)\nu(dy)$ and $Z^1$ be defined by $Z(t)=Z^1(t)+Z^2(t)$. Then $Z^1, Z^2$ are independent and $Z^2$ is a compound Poisson process. Let us consider 
the following equations
\begin{gather*}
dR(t)=F(R(t))dt+\langle G(R(t-)),dZ(t)\rangle, \quad R(0)=x,\\[1ex]
dR^1(t)=F(R^1(t))dt+\langle G(R^1(t-)),dZ^1(t)\rangle, \quad R^1(0)=x.
\end{gather*}
For the exit time $\tau_1$ of $R^1$ from the set $[(x-\varepsilon)\vee 0,x+\varepsilon]$ and the first jump time $\tau_2$ of $Z^2$ we can find
$T>0$ such that $\mathbb{P}(\tau_1>T, \tau_2<T)=\mathbb{P}(\tau_1>T)\mathbb{P}(\tau_2<T)>0$. On the set $\{\tau_1>T, \tau_2<T\}$ we have $R(\tau_2-)=R^1(\tau_2-)$ and therefore 
$$
R(\tau_2)=R^1(\tau_2-)+\langle G(R^1(\tau_2-)),\triangle Z^2(\tau_2)\rangle<-\frac{c}{2}.
$$
In the last inequality we used \eqref{pani w szpileczkach}. This contradicts the positivity of $R$. \hfill $\square$


\begin{thebibliography}{99}

\bibitem{Alfonsi1}
Alfonsi A.: Affine Diffusions and Related Processes: Simulation, Theory and Applications, (2015), Springer,




\bibitem{Bandorff-NielsenShepard}
Barndorff-Nielsen O.E., Shephard N.: Modelling by L\'evy processes for financial econometrics, (2001), 
In: Barndorff-Nielsen, O.E., et al. (eds.) L\'evy Processes: Theory and Applications,  283 - 318.
Birkhäuser, 

\bibitem{BarskiZabczykCIR}
Barski M., Zabczyk J.: On CIR equations with general factors, (2020), {\it SIAM J.Financial Mathematics}, 11,1,131-147,



\bibitem{BarskiZabczyk}
Barski M., Zabczyk J.: Bond Markets with L\'evy Factors, (2020), Cambridge University Press, 


\bibitem{BarskiZabczykArxiv}
Barski M., Zabczyk J.: A note on generalized CIR equations, (2021), {\it Communications in Information and Systems},
21, 2, 209-218,


\bibitem{CheriditoFilipovicKimmel}
Cheridito P., Filipovi\'c D., Kimmel R.L.: A note on the Dai - Singleton canonical representation of Affine Term Structure Models, (2010), {\it Mathematical Finance }, 20, 3, 509-519,


\bibitem{CuchieroFilipovicTeichmann}
Cuchiero C., Filipovi\'c  D., Teichmann J.: Affine models, (2010), {\it Encyclopedia of Quantitative Finance}, 

\bibitem{CuchieroTeichmann}
Cuchiero C., Teichmann J.: Path properties and regularity of affine processes on general state spaces, (2013),
S\'eminaire de Probabilit\' es XLV,


\bibitem{DaiSingleton}
Dai Q., Singleton K.: Specification Analysis of Affine Term Structure Models, 
(2000),  {\it The Journal of Finance}, 5, 1943-1978,


\bibitem{DuffieFilipovicSchachermeyer}
Duffie D., Filipovi\'c D., Schachermayer W.: Affine processes and applications in finance, (2003),  {\it The Annals of Applied Probability}, 13(3), 984-1053,

\bibitem{DuffieGarleanu}
Duffie, D., G\^arleanu, N.: Risk and valuation of collateralized debt obligations, (2001),  {\it Financial Analysts Journal}, 57,
41-59,


\bibitem{CIR}
Cox, I., Ingersoll, J., Ross, S.: A theory of the Term Structure of Interest Rates, (1985), {\it Econometrica}, 53,
385-408,



\bibitem{Feller}
Feller W.: An Introduction to Probability Theory and Its
Applications vol II, John Willey and Sons (1970);

\bibitem{FilipovicATS}
Filipovi\'c, D.: A general characterization of one factor affine term structure models, (2001), {\it Finance and Stochastics},
5, 3, 389-412,

\bibitem{JiaoMaScotti}
Jiao Y., Ma C., Scotti S.: Alpha-CIR model with branching processes in sovereign interest rate modeling, (2017)
 {\it Finance and Stochastics}, 21, 789-813,

\bibitem{KawazuWatanabe}
Kawazu K., Watanabe S.: Branching processes with immigration and related limit theorems, (1971)
{\it Theory Probab. Appl.}, 16, 36–54,

\bibitem{Rusinek}
Rusinek, A.: Invariant measures for forward rate HJM model with L\'evy noise. Preprint IMPAN 669
(2006),  http://www.impan.pl/Preprints/p669.pdf

\bibitem{Sato} 
Sato, K.I.: L\'evy Processes and Infinite Divisible
   Distributions, Cambridge University Press (1999),


\bibitem{Vasicek}
Vasi\v cek, O.: An equilibrium characterization of the term structure, (1997),   {\it Journal of Financial Economics}, 5, (2), 177-188.



\end{thebibliography}
\end{document}